\documentclass[12pt]{article}
\usepackage{amscd,amsfonts,amssymb,amscd,amsmath,latexsym}
\usepackage[hypertex,backref]{hyperref}
\usepackage{pslatex}
\makeatletter
\renewcommand{\subsubsection}{\@startsection
{subsubsection}
{1}
{0mm}
{0mm}
{0mm}
{\normalfont\normalsize\itshape}}
\makeatother

\newtheorem{theorem}{Theorem}[section] 
\newtheorem{prop}[theorem]{Proposition}
\newtheorem{lem}[theorem]{Lemma}
\newtheorem{ddd}[theorem]{Definition}
\newtheorem{kor}[theorem]{Corollary}

\newtheorem{axiom}[theorem]{Axiom}

\newcommand{\MSpin}{\mathbf{MSpin}}

\newcommand{\Th}{{\mathbf{Th}}}

\newcommand{\reals}{\mathbb{R}}

\newcommand{\bi}{\mathbf{i}}

\newcommand{\bl}{\mathbf{l}}

\newcommand{\naturals}{\mathbb{N}}
\newcommand{\integers}{\mathbb{Z}}

\newcommand{\into}{\hookrightarrow}

\newcommand{\iso}{\cong}

\newcommand{\forget}[1]{}

{\catcode`@=11\global\let\c@equation=\c@theorem}

\newcommand{\Z}{\mathbb{Z}}

\newcommand{\bE}{{\bf E}}

\newcommand{\proof}{{\it Proof.$\:\:\:\:$}}

\newcommand{\R}{\mathbb{R}}

\newcommand{\bL}{{\bf L}}

\newcommand{\C}{\mathbb{C}}

\newcommand{\cH}{\mathcal{H}}

\newcommand{\cK}{\mathcal{K}}

\newcommand{\Hom}{{\tt Hom}}

\newcommand{\rk}{{\tt rank}}
\newcommand{\im}{{\tt im}}

\newcommand{\bc}{\mathbf{c}}

\newcommand{\id}{{\tt id}}

\def\imath{{i}}

\def\hB{\hspace*{\fill}$\Box$ \newline\noindent}

\def\hB{\hspace*{\fill}$\Box$ \\[0cm]\noindent}

 \newcommand{\cG}{\mathcal{G}}

\newcommand{\pr}{{\tt pr}}

\newcommand{\ev}{{\tt ev}}
\newcommand{\bV}{\mathbf{V}}

\newcommand{\bh}{{\mathbf{h}}}

\title{On the topology of $T$-duality}
\author{Ulrich Bunke and Thomas Schick
\thanks{Mathematisches Institut, Universit{\"a}t G{\"o}ttingen,
Bunsenstr. 3-5, 37073 G{\"o}ttingen, GERMANY,  bunke@uni-math.gwdg.de,
schick@uni-math.gwdg.de}
}

\begin{document}


\maketitle

\begin{abstract}
We study a topological version of the $T$-duality relation between pairs consisting of
a principal
$U(1)$-bundle equipped with a degree-three integral cohomology class.
We describe the homotopy type of a classifying space for such
pairs and show that it admits a selfmap which implements a
$T$-duality transformation.

We give a simple derivation of  a $T$-duality isomorphism for certain twisted
cohomology theories. We conclude with some explicit computations
of  twisted $K$-theory groups and discuss an example of iterated
$T$-duality for higher-dimensional torus bundles.
\end{abstract}

\tableofcontents
\parskip3ex

\section{Introduction}

\subsection{Summary}

  \subsubsection{} In this paper, we describe a new approach to topological $T$-duality
  for $U(1)$-principal bundles $E\rightarrow B$ ($E$ is the background space time) equipped
  with  degree-three  cohomology
  classes $h\in H^3(E,\Z)$ (the $H$-flux in the language of the
  physical literature). 

\subsubsection{}
We first define a $T$-duality relation between such pairs
using a Thom class on an
associated $S^3$-bundle. Then 
we introduce the functor $B\mapsto P(B)$ which associates to each space  the set
of isomorphism classes of pairs. 
We construct a classifying  space $R$  of $P$ and
 characterize its homotopy type. It admits a homotopy class of selfmaps
$T:R\rightarrow R$ which implements a natural $T$-duality transformation
$P\rightarrow P$ of order two. This transformation  maps a class of pairs
$[E,h]\in P(B)$ to a canonical class $[\hat E,\hat h]\in P(B)$ of $T$-dual pairs.

 We conclude in particular that our
  definition of topological $T$-duality essentially coincides with
  previous definitions, based on integration of cohomology classes
  along the fibers.


\subsubsection{}
  
We describe an axiomatic framework for a twisted generalized cohomology
  theory $h$. We further introduce the condition of $T$-admissibility.
Examples of $T$-admissible theories are the  usual twisted de Rham cohomology and twisted K-theory.
  For a $T$-admissible generalized twisted cohomology theory $h$ we prove a $T$-duality isomorphism between
  $h(E,c)$ and $h(\hat E,\hat c)$, where $(E,c)$ and $(\hat E,\hat
  c)$
are $T$-dual pairs.
  
\subsubsection{}
  We compute a number of examples. Iterating the construction of
  $T$-dual pairs, we can
  define duals of certain higher dimensional torus bundles. We show
  that with our definition of duality the isomorphism type of the
  dual of a torus bundle, even if it exists, is not always uniquely determined.

\subsubsection{}
We thank the referees for their useful comments, in particular with respect
to the presentation and the physical interpretation of our results.

\subsection{Description of the results}

\subsubsection{}
In this paper we try to explain our understanding of the results of the recent paper
\cite{bem} and parts of \cite{bhm} and \cite{mr} (Sec. 4.1) by means
of elementary algebraic topology. The notion of $T$-duality originated in string
theory. Instead of providing an elaborate historical account of
$T$-duality here we refer to the two papers above and the literature cited
therein. In fact, the first paper which studies $T$-duality  
is in some sense \cite{rr}. We will explain the relation with the
present paper later in this introduction.

\subsubsection{}
However, a few motivating words what this paper is about, and more
importantly what it is not about, are in order.

$T$-duality first came up in physics in the following situation. The space
$E$ appears as a part of a ``background space-time". The
cohomology class $h\in H^3(E,\Z)$ describes the Flux for a Neveu-Schwarz
$3$-form gauge potential $H$. In connection with $T$-duality, the case where $E$ admits a free
$U(1)$-action and thus has the structure of a principal
$U(1)$-bundle, is of particular interest. The natural generalization
is a space with a free action of a higher-dimensional torus.
Then $E$ is a
$U(1)^k$-principal bundle. In such a situation, for physical reasons,
one expects to find a dual bundle with a dual flux (i.e.~cohomology class) roughly by replacing each fiber by
the dual fiber, the so called $T$-dual. The $T$-dual should share
many properties with the original bundle. In particular 
one expects that certain  twisted
cohomology groups are isomorphic.

In the physical situation the spaces
come with geometry. When passing to the dual, the metric on the
fiber should be replaced by the dual metric on the dual fiber. A lot
of the literature about $T$-duality and its relation to mirror
symmetry have the geometry as a major ingredient, and they focus on
situations in which the dimension of the fiber and the base
coincide. One of the basic contributions in this context is \cite{syz}.

\subsubsection{}

In the present paper we will completely disregard the geometry and
metrics. This also explains the title  ``\emph{topological}
$T$-duality''. We are only interested in the resulting topological
type. Moreover, we adopt a mathematical definition of the $T$-duality
relation by
simply declaring certain cohomological properties which are 
expected for physical reasons. 

This approach works best for $U(1)$-bundles. So we will concentrate on
those for most of the paper with the exception of Section
\ref{sec:iterated-t-duality}, where we study torus bundles by considering them as
iterated $U(1)$-bundles.

\subsubsection{}

In the present paper we study $T$-duality for principal $U(1)$-bundles
equipped with an integral cohomology class of degree $3$.
We will call such data a \emph{pair} (Definition \ref{pdef}).
We first introduce $T$-duality as a \emph{relation} between  pairs
(Definition \ref{dude}) (in particular, a given pair can have several
$T$-dual pairs). The paper \cite{bem} works almost exactly in
the same setting: it also starts with a pair and defines what a {dual
pair} is (via a construction which involves some choices, so again is
not unique). This
definition unfortunately is not very precise, since torsion in the
  cohomology is neglected. In Section \ref{sec:an-example-where} we
  show by an example that it is necessary to take the torsion into account if
  one studies e.g. the $T$-duality isomorphism for
  twisted K-theory.

\subsubsection{} 

At a first
glance our definition of $T$-duality, which is based on a Thom class on an
auxiliary $3$-sphere bundle, looks quite
different from the definition given in \cite{bem}, which relied  on
integration over the fiber. The link between the two definitions is provided by an explicit
universal example over a universal base space $R$ for our definition of
$T$-duality. Using some non-trivial
calculations in this universal example we will obtain
a complete characterization of the $T$-dual (according to our
definition) by topological
invariants, which contains in particular the same kind of integration
over the fibre as in
the older notions of $T$-duality. We will show that, up to passing to real cohomology,
the $T$-duality of \cite{bem} is characterized by the same topological
invariants. Therefore, we
can eventually conclude that our definition is essentially equivalent to the one
used there (see \ref{zu71}, \ref{zu72}).

Later in the present paper we will understand $T$-duality as a map which
associates to an isomorphism class of  pairs a canonical dual
isomorphism class of pairs  in a two-periodic manner.
This in particular reproves the result of \cite{bem} that each pair
admits a $T$-dual.

\subsubsection{}

A third definition of $T$-duality is given in
 \cite{rr}
(compare also \cite{mr}, 4.1) or in \cite{bem}. In \cite{rr}, the main object is a
continuous trace algebra $A$ with an $\R$-action such that its
spectrum $X(A)$ is
a free $U(1)\cong \R/\Z$-space. To $A$ we can associate a
pair $(X(A),h(A))$ consisting of the $U(1)$-bundle
$X(A)\rightarrow X(A)/U(1)$ and the Dixmier-Douady class
$h(A)\in H^3(X(A),\Z)$. Vice versa each pair can be realized in this
way. With an appropriate notion of Morita equivalence we have a
bijection of equivalence classes of such algebras and isomorphism
classes of pairs.

In \cite{rr} it is shown that the cross product $\hat A:=A\rtimes \R$
is again a continuous trace algebra with $\R$-action (the latter $\R$ is in
fact the dual group of $\R$) of the same type as above.
It follows from the comparison of the topological invariants
of the pairs $(X(\hat A),h(\hat A))$ and the dual pair
$(\widehat{X(A)},\widehat{h(A)})$
and the naturality of the constructions with respect to the change of
the base spaces
 that our  notion of $T$-duality of pairs
indeed corresponds to the cross product in \cite{rr}. 

It is well known  that $\widehat{\hat A}$ is Morita equivalent to
$A$. This fact is reflected in our picture by the result that
$T$-duality is two-periodic. 

\subsubsection{}

Given a base space $B$, we study the set $P(B)$ of isomorphism classes of
pairs $(E,h)$ over $B$, where $E\rightarrow B$ is  of a
$U(1)$-principal bundle and $h$ a class $h\in
H^3(E,\Z)$. It turns out that the contravariant set-valued
functor $B\mapsto P(B)$ can be represented by a space $R$, the
classifying space of pairs. The $T$-duality can then be considered as
a natural transformation $T:P\rightarrow P$ of functors, and it is
represented by a homotopy class of maps $T:R\rightarrow R$.

\subsubsection{}

Our first main result (Theorem \ref{main1}) is the characterization of
 the homotopy type of $R$ as the  homotopy
fibration
$$K(\Z,3)\rightarrow R\rightarrow K(\Z,2)\times K(\Z,2)$$
which is classified by
$$\pr_1^* c\cup\pr_2^* c\in H^4(K(\Z,2)\times K(\Z,2),\Z).$$
Here $K(\Z,n)$ is the Eilenberg-MacLane space, i.e.~characterized by
the property that $\pi_k(K(\Z,n))=0$ if $k\ne n$ and
$\pi_n(K(\Z,n))=\Z=H^n(K(\Z,n),\Z)$. In particular, we can choose
$K(\Z,2)=\C P^\infty$. The class $c\in H^2(K(\Z,2),\Z)$ is the
canonical generator. How fibrations are classified is recalled in
\ref{sec:universal_fibrations}. 

\subsubsection{}

The space $R$ carries a universal pair, and the map $T$ will of course
have the property to represent the universal dual pair (Definition
\ref{dpa}).

The classifying space $R$ in fact already  appears  in \cite{rr} (proof of
Theorem. 4.12). It is used there in order to simplify the
verification of the relation of topological invariants which
corresponds to the assertion of Lemma
\ref{lem2}.

\subsubsection{}

As observed in many places, $T$-duality comes with isomorphisms in
certain twisted generalized cohomology theories. In fact, the
calculation of such twisted generalized cohomology groups in terms of
the (perhaps easier to understand) generalized cohomology groups of
the dual is one (topological) motivation for the study of $T$-duality. If $(E,h)$ and $(\hat E,\hat h)$
are pairs over $M$, and in particular $E$ and $\hat E$ are principal
$U(1)$-bundles over $M$, which are dual to each other, than
(as shown e.g.~in \cite{bem})
there is an isomorphism (of degree $-1$) in twisted complex $K$-theory
$K_*(E,h)\cong K_{*-1}(\hat E,\hat h)$ or of  real twisted cohomology
 $H(E,\R,h)\cong H(\hat E,\R,\hat h)$.
These isomorphisms are implemented by explicit $T$-duality
transformations (Definition \ref{trd})
which are constructed out of the diagram
\begin{equation}\begin{array}{ccccc}
&&E\times_B\hat E&&\\
&p\swarrow&&\hat p\searrow&\\
E&&q\downarrow&&\hat E\\
&\pi \searrow&&\hat \pi \swarrow\\
&&B&&
\end{array}
\ \end{equation}
using standard operations in twisted cohomology (like pull-back and
integration over the fiber)\footnote{In the $C^*$-algebraic context of
\cite{rr}, \cite{mr} the $T$-duality isomorphism is given by Connes'
Thom isomorphism for crossed products with $\R$.
}.

\subsubsection{}
We say that a twisted generalized cohomology theory is $T$-admissible
if the $T$-duality transformation is an isomorphism in the special
case of the pair $(U(1)\rightarrow *,0)$.
Our second main result is the observation (Theorem \ref{main2})
that the  $T$-duality transformation for a $T$-admissible twisted
generalized cohomology theory is an isomorphism, and that this fact is an
easy consequence of the Mayer-Vietoris principle.

\subsubsection{}
In order to produce a precise statement we fix the axioms
for a twisted generalized  cohomology theory in Subsection \ref{fff}.
In doing so we add some precision to the statements in \cite{bem},
in particular to the observation that the  Chern character preserves the
$T$-duality transformation (\cite{bem}, 1.14).

The main point
  is that the cohomology class $h\in H^3(E,\Z)$  only determines the
  \emph{isomorphism class} of a twist and so the isomorphism class of
  $K(E,h)$ or $H(E,\R,h)$ as an abstract group. In order to be able to
  say that
  the Chern character is a transformation between twisted cohomology
  theories one must use the same explicit objects to twist $K$-theory
  as one uses to twist real cohomology. In order to twist complex $K$-theory one
  usually  considers a  principal
  $PU$-bundle  (but not a
  three-form as in \cite{bem}). More details on twisted K-theory can
  be found in  \cite{atiyahsegal}. On the other hand, three-forms are
  usually used to twist real (de Rham) cohomology.
We do not know any natural way to relate these two kinds of twists
(but look at \cite{fht1}, proof of Prop. 3.5, which perhaps solves
this problem).
In a previous paper \cite{bunkeschick041} we have constructed
versions of twisted $K$-theory and twisted real cohomology where the
twists in both cases are Hitchin gerbes. For these versions of twisted
cohomology theories the Chern character is indeed a
natural transformation and preserved by $T$-duality. Since this gives a framework
to work simultaneously with twisted K-theory and twisted cohomology,
we propose to use Hitchin gerbes in this context. In the paper,
however, we simply assume that the twists $\cH$ and the twisted generalized
cohomology theory $h$ satisfies certain natural axioms, and then we go
on to prove a natural
$T$-duality isomorphism $h(E,\cH)\xrightarrow{\iso} h(\hat E,\hat \cH)$ for
any theory which satisfies these axioms and for dual pairs $(E,\cH)$
and $(\hat E,\hat \cH)$. 

\subsubsection{}

For the purpose of illustration
we perform some calculations of twisted $K$-theory.
For three-manifolds we obtain a complete answer in Subsection
\ref{dre} (compare with the partial results of
\cite{mi}). We demonstrate the $T$-duality isomorphism in twisted
$K$-theory for $U(1)$-principal bundles over surfaces by explicit
calculation.

\subsubsection{}

It is a natural question if $T$-duality can be generalized to
 principal $U(1)^k$-bundles for $k>1$. As observed in \cite{bhm}
and \cite{mr} not every $U(1)^k$-principal bundle has a $T$-dual in the classical
sense. {Note the remarkable observation in \cite{mr}
  Theorem 4.4.2, that in general the $T$-dual of a $U(1)^2$-principal
  bundle
equipped with a three-dimensional integral cohomology class
is a bundle of non-commutative tori.}
In the present paper we discuss the approach of defining a $T$-dual of
a higher-dimensional principal torus bundle as an iterated $T$-dual of
$U(1)$-principal bundles. We demonstrate by an example that this
approach does not lead to a unique result.

\subsubsection{}
A $U(1)$-principal bundle $E\rightarrow B$ is essentially
the same object as the free $U(1)$-space $E$. In a continuation \cite{bunkeschick043} of the
present paper we discuss a generalization of $T$-duality to the case
of $U(1)$-spaces where $U(1)$ acts with at most finite
stabilizers. For applications to physics, this seems to be of
relevance.

\section{The classifying space of pairs}
\subsection{Pairs and the classifying space}
\subsubsection{}
Let $B$ be a topological space.
\begin{ddd}\label{pdef}
A pair $(E,h)$ over $B$ consists of a
$U(1)$-principal bundle $\pi:E\rightarrow B$
and a class $h\in H^3(E,\Z)$.
\end{ddd}

\subsubsection{}
If $f:A\rightarrow B$ is a continuous map,
then we can form the functorial
pull-back $f^*(E,h)=(f^*E,F^*h)$,
where $F$ is defined by the pull-back
 $$\begin{array}{ccc}f^*E&\stackrel{F}{\rightarrow}&E\\
\downarrow&&\downarrow\\
A&\stackrel{f}{\rightarrow}&B\end{array}\ .$$

\subsubsection{}
We say that two pairs are isomorphic (written as
$(E_0,h_0)\cong (E_1,h_1)$)
if there
exists an isomorphism of $U(1)$-principal bundles
 $$\begin{array}{ccc}E_0&\stackrel{F}{\rightarrow}&E_1\\
\downarrow&&\downarrow\\
B&=&B\end{array}$$
such that $F^*h_1=h_0$.

\subsubsection{}
Let $(E_i,h_i)$, $i=0,1$ be pairs over $B$.
We say that they are homotopic (written as
$(E_0,h_0)\sim (E_1,h_1)$), if there exists a pair
$(\tilde E,\tilde h)$ over $I\times B$ such that
$f_i^*(\tilde E,\tilde h)= (E_i,h_i)$, $i=0,1$,
where $f_i:B\rightarrow I\times B$ is given by
$b\mapsto (i,b)$. Note that we insist here in equality, it is not
sufficient for later purposes to only have an isomorphism.

\subsubsection{}

\begin{lem}\label{hom}
On pairs, the relations ``homotopy equivalence'' $\sim$ and
``isomorphism'' $\cong$ coincide.
\end{lem}
\proof
Let $(E_0,h_0)$ and $(E_1,h_1)$ be homotopic via $(\overline
E,\overline h)$. Then there is an isomorphism $\overline E\to
E_0\times [0,1]$. Using this, we immediately get an isomorphism
$F:f_0^*(\overline E,\overline h)\cong f_1^*(\overline E,\overline h)$.

Conversely, if $(E_0,h_0)$ and $(E_1,h_1)$ are isomorphic via an
isomorphism $F$, we construct the homotopy $\overline{E}:= E_0\times [0,1/2]
\cup_{F\times\id_{\{1/2\}}} E_1\times [1/2,1]$, with $\overline h$ obtained
(uniquely) using the Mayer-Vietoris sequence for the cohomology of
$\overline E$. We take the freedom to use \emph{canonical} isomorphisms
between $E_k\times \{k\}$ and $E_k$, $k=0,1$.
\hB

\subsubsection{}

\begin{ddd}
  By $K(\Z,n)$ we denote the Eilenberg-Mac Lane space characterized
  (upto homotopy equivalence) by its homotopy groups $\pi_k(K(\Z,n))=0$
  if $k\ne n$, $\pi_n(K(\Z,n))=\Z$. Recall that for an arbitrary space
  $X$ the cohomology with $\Z$-coefficients $H^n(X,\Z)$ can be identified with homotopy classes of maps from
  $X$ to $K(\Z,n)$ (denoted by $[X,K(\Z,n)]$), a fact we are going to
  use frequently.
\end{ddd}

As a model for $K(\Z,1)$ we choose $U(1)$.
As a model for $K(\Z,2)$ we can
choose $\C P^\infty$.
Let $q:
U\rightarrow K(\Z,2)$ be the universal $U(1)$-principal bundle. If we
choose $K(\Z,2)=\C P^\infty$, we can choose $U:=S(\C^\infty)$, the unit
sphere in $\C^\infty=\bigcup_{n\in\naturals} \C^n$, and $p$ factors
out the canonical $U(1)$-action on $\C^\infty$.
Furthermore, let $LK(\Z,3)$ be the free loop space of
$K(\Z,3)$. This space admits an action of $U(1)$ by
$u\gamma(t):=\gamma(u^{-1}t)$
for $\gamma\in LK(\Z,3)$ and $u,t\in U(1)$.

\begin{ddd}
We define the space $R$ as the total space of the associated bundle
$$\bc\colon (R:=U\times_{U(1)}LK(\Z,3))\rightarrow K(\Z,2)\ .$$
\end{ddd}
Note that $R$ is well defined up to homotopy equivalence. We
consider $\bc$ also as a cohomology class $\bc\in H^2(R,\Z)$.

\subsubsection{}
Over $R$ we have the $U(1)$-principal bundle $\pi\colon (\bE:=\bc^*U)\rightarrow
R$ with first Chern class $\bc\in H^2(R,\Z)$. Furthermore, we have a canonical map
\begin{equation*}
\bh:\bE\rightarrow K(\Z,3);\quad
\bh(u,[v,\gamma]):=\gamma(t),
\end{equation*}
where
$u,v\in U$, $\gamma\in LK(\Z,3)$ and $t\in U(1)$ satisfy $q(u)=q(v)=\bc([v,\gamma])$,
and $tv=u$.
Note that this is well-defined, independent of the choice
of the representative of the class $[v,\gamma]\in R$.
We consider this map also as a cohomology class $\bh\in H^3(E,\Z)$.
In this way we get a pair $(\bE,\bh)$ over $R$.

\begin{ddd}
We call this pair $(\bE,\bh)$ the universal pair.
\end{ddd}

\subsubsection{}\label{tr1}
We define the contravariant functor $P$ from the category of
topological spaces to the category of
sets which associates to the space $B$ the set $P(B)$ of isomorphism
classes of pairs
and to the map $f\colon A\rightarrow B$ the pull-back $f^*\colon P(B)\rightarrow
P(A)$.

\begin{prop}\label{cla}
The space $R$ is a classifying space for $P$. In fact, we have an
isomorphism of functors
$\Psi_{\dots }:[\dots,R]\rightarrow P(\dots)$
given by $\Psi_B([f]):=[f^*(\bE,\bh)]$ for each homotopy class of maps
$[f]\in [B,R]$ and each CW-complex $B$.
\end{prop}
\proof
It follows  immediately from Lemma \ref{hom}
that the functor $P$ is homotopy invariant.
Therefore $\Psi_{\dots}$ is a well-defined natural transformation.

Let $[E,h]\in P(B)$ be given.
Up to isomorphism, we can assume that we have a pull-back diagram of $U(1)$-principal bundles
$$\begin{array}{ccc}
E&\stackrel{C}{\rightarrow}&U\\
\downarrow&&\downarrow\\
B&\stackrel{c}{\rightarrow}&K(\Z,2)\end{array}\ .$$
We represent the class $h$ by a map
$h:E\rightarrow K(\Z,3)$.
We construct a lift $f:B\rightarrow R$ of $c$ as follows.
For $b\in B$ choose $e\in E_b$. Then we set
\begin{equation*}
f(b):=[C(e),\gamma]\in R\quad\text{with }\gamma(t)=h(te) \;\forall t\in
U(1).
\end{equation*}
 Observe that
$f(b)$ is independent of the choice of $e$. If $F\colon E\to \bE$ is
the $U(1)$-bundle map covering  $f$, then $F(e) =
(C(e),[C(e),\gamma])$ with $e$ and $\gamma$ as above.
Therefore, $\bh\circ F=h$ and we have $f^*(\bE,\bh)\cong (E,h)$.
This shows that $\Psi_B$ is surjective.

Let now $\Psi_B([f_0])=\Psi_B([f_1])$. Using Lemma \ref{hom}, we choose a
homotopy $\overline E$ over $B\times [0,1]$ between $f_0^*(\bE,\bh)$
and $f_1^*(\bE,\bh)$. The construction used for the surjectivity part
provides us with a map $\overline f\colon B\times [0,1]\to R$ such that
$\overline f^*(\bE,\bh)= (\overline E,\overline h)$. To achieve equality, we have to
choose $\overline E$ in such a way that $\overline E=\overline c^*U$
for an appropriate map $\overline c\colon B\times [0,1]\to K(\Z,2)$
(without changing the bundle at the boundary, i.e.~such that
$\overline c_k=\bc\circ f_k$). This is possible since $K(\Z,2)$ is a
classifying space for principal $U(1)$-bundles.

The construction has the property that $\overline f_k=f_k$, therefore
$\overline f$ is a homotopy between $f_0$ and $f_1$, proving that
$\Psi_B$ is injective.

\hB

\subsection{Duality of pairs}

\subsubsection{}\label{dua}
Let $\pi:E\rightarrow B$ and $\hat \pi:\hat E\rightarrow B$ be two
$U(1)$-principal bundles.
Let
\begin{equation*}
\pi: (L:=E\times_{U(1)}\C)\rightarrow B
\end{equation*}
and $\hat \pi :(\hat L:=\hat
E\times_{U(1)}\C)\rightarrow B$ be the associated complex hermitian line bundles.
We can consider $E$ and $\hat E$ as unit sphere bundles in $L$ and
$\hat L$. We form the complex vector bundle
$r:(V:=(L\oplus \hat L))\rightarrow B$ and let $r:S(V)\rightarrow B$ be
the unit sphere bundle, the fibers consisting of $3$-dimensional
spheres. $V$ being a complex vector bundle, the map $r$ is oriented.
In particular, we have an integration map
$r_!:H^3(S(V),\Z)\rightarrow H^0(B,\Z)$ (in de Rham cohomology the
corresponding map is really given by integration over the fiber).
Let $1_B$ denote the unit in the ring $H(B,\Z)$.
\begin{ddd}
A Thom class for $S(V)$ is a class
$\Th\in H^3(S(V),\Z)$ such that
$r_!(\Th)=1_B$.
\end{ddd}
If $S(V)$ admits a Thom class, then by the Leray-Hirsch theorem
its cohomology is a free $H(B,\Z)$-module generated by $1_{S(V)}$ and
$\Th$. Thom classes in general are not unique. In fact, $\Th^\prime$ is a
second Thom class if and only if $\Th-\Th^\prime=p^*d$ for some $d\in
H^3(B,\Z)$.

\subsubsection{}
Let $c,\hat c\in H^2(B,\Z)$ denote the Chern classes of $E$ and $\hat
E$. The product $\chi(V):=c\cup \hat c \in H^4(B,\Z)$ is the Euler class
of $V$.

\begin{lem}\label{eul}
The bundle $S(V)$ admits a Thom class if and only if $\chi(V)=0$.
\end{lem}
\proof
This follows from the Gysin sequence for $S(V)$. For this question the
important segment is
$$\rightarrow H^3(B,\Z)\stackrel{r^*}{\rightarrow}
H^3(S(V),\Z)\stackrel{r_!}{\rightarrow}H^0(B,\Z)\stackrel{\chi(V)}{\rightarrow}
H^4(B,\Z)\rightarrow\ .$$
\hB

\subsubsection{}
We now consider two pairs $(E,h)$ and $(\hat E,\hat h)$.
Let $i:E\rightarrow S(V)$ and
$\hat i:\hat E\rightarrow S(V)$ denote the inclusions of the
$S^1$-bundles into the $S^3$-bundle.
\begin{ddd}\label{dude}
We say that $(E,h)$ and $(\hat E,\hat h)$ are dual to each other if there exists a
Thom class $\Th$ for $S(V)$ such that
$h=i^*\Th$ and $\hat h=\hat i^*\Th$.
\end{ddd}

\subsubsection{}

Let $\pi:E\rightarrow B$ and $\hat \pi:\hat E\rightarrow B$ be given $U(1)$-principal bundles with first
Chern classes $c$ and $\hat c$.
Then \ref{eul} has the following consequence.
\begin{kor}\label{var}
There exists $h\in H^3(E,\integers)$ and $\hat h\in H^3(\hat E,\Z)$
such that $(E,h)$ and $(\hat E,\hat h)$ is a dual pair, if and only if
$c\cup \hat c=0$. If such a dual pair exist, then any other
has the form $(E,h+\pi^*b)$ and $(\hat E,\hat h+\hat \pi^*b)$ for some
$b\in H^3(B,\Z)$.

\end{kor}

\subsubsection{}
Let $(E,h)$ and $(\hat E,\hat h)$ be dual pairs.
We consider the following part of the Gysin sequence for $E$
$$\rightarrow H^1(B,\Z)\stackrel{c}{\rightarrow}
H^3(B,\Z)\stackrel{\pi^*}{\rightarrow} H^3(E,\Z)\rightarrow\ .$$
We observe the following consequence of \ref{var}.
\begin{kor}\label{corol:difference_between_different_duals}
If $(E,h)$ is dual to $(\hat E,\hat h)$ and also to $(\hat E,\hat h')$, then
we have
$\hat h'-\hat h=\hat \pi^*(c\cup a)$ for some $a\in H^1(B,\Z)$.
\end{kor}

\begin{lem}\label{cde}
If $(E,h)$ is dual to $(\hat E,\hat h)$,
then $c=-\hat \pi_!(\hat h)$ and $\hat c=-\pi_!(h)$.
\end{lem}
\proof
We defer the proof to \ref{lem2}. It follows from the calculation of
the cohomology in the universal situation.
\hB

\begin{lem}
  Let $(E,h)$ be dual to $(\hat E,\hat h)$. Consider the fiber product
\begin{equation}\label{cdnew}\begin{array}{ccccc}
&&E\times_B\hat E&&\\
&p\swarrow&&\hat p\searrow&\\
E&&q\downarrow&&\hat E\\
&\pi \searrow&&\hat \pi \swarrow\\
&&B&&
\end{array}
\ .\end{equation}
Then $p^*h=\hat p^*\hat h$.  
\end{lem}
\proof This is the parameterized version of the situation considered
later in \ref{poi}. In particular, we have a homotopy
$h:I\times E\times_M\hat E\rightarrow S(V)$ from $i\circ p$ to $\hat
i\circ \hat p$, where $i\colon E\to S(V)$ and $\hat{i}\colon \hat E\to
S(V)$ are the canonical inclusions into the sphere bundle of the complex vector bundle $V$
associated to $E$ and $\hat E$, then $p^*h=p^*i^*\Th=\hat p^*\hat
i^*\Th=\hat p^*\hat h$.
\hB

\subsubsection{}\label{zu71}
We are now in the situation to compare our definition of $T$-duality
with the definition used in \cite{bem}, Section 3.1. When interpreted in
cohomological terms instead of using the language of differential forms,
\cite{bem} constructs to a given pair
$(E,h_\R)$ (where $h_\R\in \im(H^3(E,\Z)\to H^3(E,\R))$ is an real
cohomology class with integral periods) another pair $(\hat E,\hat h_\R)$,
again with $\hat h_\R\in H^3(\hat E,\R)$.

Let $c$ be the first Chern class of $E$ and use the notation of \eqref{cdnew}. 
By $c_\R$ we denote the
image of $c$ in $H^2(B,\R)$.

The
construction in \cite{bem} depends on a few choices, in particular the choices of
connections. An integral lift $h\in H^3(E,\Z)$ of
$h_\R$ uniquely determines the isomorphism class of the
$U(1)$-principal bundle $\hat E$ with Chern class $\hat c:=\pi_!(h)$. The cohomology class $\hat h_\R$ is
then 
determined up to addition of a class of the form $\hat \pi^*(c_\R\cup b)$
with some $b\in H^1(B,\R)$.

In \cite{bem}, 3.1 it is shown that
$\pi_*(h_\R)=\hat c_\R$ and $\hat\pi_*(\hat h_\R)=c_\R$. 
These formulas differ from those of Lemma \ref{cde} by some signs. The
reason is that in \cite{bem} the dual bundle is considered with the
opposite
$U(1)$-action. 
In \cite{bem} it
 is also shown that $p^*h_\R=\hat
p^*\hat h_\R$.

We will now prove that up to addition of classes of the form $\hat
\pi^*(c_\R\cup b)$ for $b\in H^1(B,\R)$ the class $\hat h_\R\in H^3(\hat E,\R)$ is
uniquely determined by these properties. Since our $T$-duality pairs
share these properties, we conclude that (upon passing to real
cohomology) they are dual in the sense of \cite{bem}. It then follows
also that $\hat h_\R$ can be chosen with integral periods and with an
integral lift $\hat h$ such that $\hat \pi_*\hat h=c$, since we
construct an integral lift of some representative. This assertion is also
implicit in \cite{bem}, but without a detailed proof. Note also that
the ambiguity in the dual class $\hat h$ is exactly parallel to the
ambiguity in the construction of \cite{bem}.

\subsubsection{}\label{zu72}

To prove that $\hat h_\R$ is determined by the properties $\hat
\pi_*\hat h_\R=c_\R$ and $\hat p^*\hat h_\R=p^*h_\R$ we consider the
following web of Gysin sequences
 for the $U(1)$-principal bundles
$p$, $\hat p$, $\pi$ and $\hat \pi$. Every row and every column is
exact, and by the naturality of the Gysin sequence every square
commutes. We use cohomology with real coefficients throughout, but the
diagram is of course also correct with integral coefficients.

\begin{equation}\label{eq:bigdg}
  \begin{CD}
    0@>>> H^1(B) @>{\hat\pi^*}>> H^1(\hat E) @>{\hat \pi_!}>> H^0(B) @>{\cup \hat
      c}>> H^2(B)\\
    @VVV @VV{\cup c}V @VV{\cup\hat\pi^*c}V @VV{\cup c}V @V{\cup c}VV\\
   H^1(B) @>{\cup \hat{c}}>>  H^3(B) @>{\hat\pi^*}>> H^3(\hat E)
   @>{\hat \pi_!}>> H^2(B) @>{\cup\hat c}>> H^4(B)\\
    @VVV @VV{\pi^*}V @VV{\hat p^*}V @VV{\pi^*}V @VV{\pi^*}V\\
   H^1(E) @>{\cup c}>> H^3(E) @>{p^*}>> H^3(E\times_B\hat E)
      @>{p_!}>> H^2(E)@>{\cup c}>> H^4(E)
  \end{CD}
\end{equation}

Assume that $\hat h,\hat h'\in H^3(\hat E)$ both satisfy the above
equations, and set $d:=\hat h-\hat h'$. It follows that $\hat \pi_! d=0\in
H^2(B)$, and that $\hat p^*d=0$. The second property implies that
there is a lift $l\in H^1(\hat E)$ with $d=l\cup \hat \pi^*c$. Set
$n:=\hat \pi_! l\in H^0(B)$. Without loss of generality we can assume that
$B$ is connected (else we work one component at a time). Now, only two
possibilities remain (since \cite{bem} uses real coefficients, where no
torsion phenomena occur).
\begin{enumerate}
\item Either $n=0$, then $l=\hat \pi^*a$ for a suitable $a\in H^1(B)$,
  and consequently $\hat h-\hat h'=d = \hat \pi^*(c\cup a)$, which is
  exactly what we want to prove.
\item If $n\ne 0$, then $c_\R=0$, since $nc_\R=\hat \pi_! d=0$. In this
  case, $\hat \pi^*c_\R=0$ and therefore also $\hat h-\hat h'=d=l\cup
  \hat\pi^*c_\R =0$.
\end{enumerate}

\subsubsection{}
Let us fix $(E,h)$.
\begin{theorem}\label{uni}
The equivalence class of pairs which are dual to $(E,h)$ is uniquely
determined.
\end{theorem}
\proof
By Lemma \ref{cde} the isomorphism class of the underlying $U(1)$-bundle $\hat E$ of a
pair dual to $(E,h)$ is determined by the first Chern class $\hat
c:=\pi_!(h)$. If $(\hat E,\hat h)$ and $(\hat E,\hat
h^\prime)$ are both dual to $(E,h)$, then by Corollary
\ref{corol:difference_between_different_duals} $\hat h^\prime-\hat h=\hat
\pi^*(c\cup a)$ for some $a\in H^1(B,\Z)$.
It remains to show that there exists an automorphism of
$U(1)$-principal bundles
$$\begin{array}{ccc}\hat E&\stackrel{U}{\rightarrow}&\hat E\\
\downarrow&&\downarrow\\
B&=&B\end{array}$$
such that $U^*\hat h=\hat h^\prime$.
Any automorphism $U$ is given by multiplication by a suitable $g:B\rightarrow U(1)$.
Then we can factor $U$ as the composition
$$\hat E\stackrel{(\hat \pi,\id)}{\rightarrow} B\times \hat E\stackrel{g\times \id}{\rightarrow}U(1)\times
\hat E\stackrel{m}{\rightarrow} \hat E\ ,$$
where $m$ is given by the principal bundle structure.
Observe that we have the pull-back diagram
$$\begin{array}{ccc}U(1)\times \hat E&\stackrel{m}{\rightarrow}&\hat E\\
\pr_2\downarrow&&\pi\downarrow\\
\hat E&\stackrel{\pi}{\rightarrow} &B\end{array}\ .$$
Using $H^3(U(1)\times \hat E)=\pr_2^* H^3(\hat E)\oplus o_{U(1)}\times
\pr_2^*H^2(\hat E)$ where $o_{U(1)}\in H^1(U(a))$ is the canonical generator, naturality of integration over the fiber, and the
split of $\pr_2$, we obtain
$$m^*(\hat h)=\pr_2^*(\hat h)\oplus o_{U(1)}\times \hat \pi^*\hat \pi_!(\hat h)\ .$$
Note that $[B,U(1)]\cong H^1(B,\Z)$ via $[g]\mapsto
g^*o_{U(1)}=:a(g)\in H^1(B,\Z)$.

Now we return to the construction of $U$ (and therefore $g$) with
$U^*\hat h=\hat h'$. To achieve this, choose $g$ corresponding to $a\in
H^1(B,\Z)$ such that $\hat h-\hat h'=\hat\pi(c\cup a)$. Using $\hat \pi_!(\hat h)=-c$ we get
$(g\times \id)^*m^*(\hat h)=-a\times \pi^*c+\pr_2^*(\hat h)$.
Finally  $U^*(\hat h)=\hat h-\pi^*(c\cup a)$. \hB

\subsection{The topology of $R$}

\subsubsection{}
It is a topological fact that the universal bundle with fiber $K(\Z,3)$
is
$$K(\Z,3)\rightarrow PK(\Z,4)\rightarrow K(\Z,4)\ ,$$
where $PK(\Z,4)$ is the path space of $K(\Z,4)$,
i.e. the space of all path in $K(\Z,4)$ starting in the base
point. The map to $K(\Z,4)$ is given by evaluation at the end point.
The fiber of this evaluation over the base point
is the based loop space $\Omega K(\Z,4)$,
which serves here as a model for the homotopy type $K(\Z,3)$.

\subsubsection{}\label{sec:universal_fibrations}
If $B$ is a space, then bundles over $B$ with (``oriented'') fiber $K(\Z,3)$ are
classified
by homotopy classes of maps $[B,K(\Z,4)]$, i.e. by cohomology classes
in $H^4(B,\Z)$. The homotopy type of such a bundle is determined by
such maps up to self homotopy equivalences of $B$ and of $K(\Z,4)$,
i.e.~upto the action of self homotopy equivalences of $B$ and up to
multiplication by $-1$ on $H^4(B;\Z)$.

We consider a bundle
$$K(\Z,3)\rightarrow F\rightarrow B$$
which is classified by $\kappa\in H^4(B,\Z)$.
For simplicity we assume that $B$ is connected and simply connected.
Then $\kappa$ can be read off from the differential
$d_4^{0,3}$ in the Serre spectral sequence for the bundle.
By the Hurewicz theorem, the relevant part of the $E_4$-page looks like
$$\begin{array}{|c|c|c|c|c|c|}\hline
3&\Z&0&H^2(B,\Z)&H^3(B,\Z)&H^4(B,\Z)\\\hline
2&0&0&0&0&0\\\hline
1&0&0&0&0&\\\hline
0&\Z&0&H^2(B,\Z)&H^3(B,\Z)&H^4(B,\Z)\\\hline
X&0&1&2&3&4\\\hline
\end{array}\ .$$
The differential $d_4^{0,3}:\Z\rightarrow H^4(B,\Z)$ is  multiplication with $\kappa$.

\subsubsection{}
The main result of the present subsection is the determination of the
homotopy type of $R$. Let $z\in H^2(K(\Z,2),\Z)$ be the canonical generator.
By the K{\"u}nneth theorem, the cohomology of $K(\Z,2)\times K(\Z,2)$ is the polynomial ring in
two generators $c=\pr_1^*z$ and $\hat c:=\pr_2^*z$, i.e.
$H(K(\Z,2)\times K(\Z,2),\Z)=\Z[c,\hat c]$.
\begin{theorem}\label{main1}
$R$ is the total space of a bundle
\begin{equation}\label{eq3}K(\Z,3)\rightarrow R\rightarrow
  K(\Z,2)\times K(\Z,2)\end{equation}
which is classified by
$c\cup \hat c\in H^4(K(\Z,2)\times K(\Z,2),\Z)$.
\end{theorem}

\subsubsection{}
To prove Theorem \ref{main1}, we first compute the homotopy groups $\pi_i(R)$.
Observe that $S^0$ and $S^1$ admit only one isomorphism class of
pairs. This implies that $R$ is connected and simply connected.
This observation also frees us from basepoint considerations.


\begin{lem}
The homotopy groups of $R$ are given by
$$\pi_i(R)=\left\{\begin{array}{cc}0&i\not\in\{2,3\}\\
\Z\oplus \Z&i=2\\
\Z&i=3\end{array}\right.\ .$$\end{lem}
\proof
We first observe that there is exactly one isomorphism class of pairs
over $S^i$ for $i\ge 4$,
namely $(U(1)\times S^i\rightarrow S^i,0)$. This implies that
$\pi_i(R)=0$ for $i\ge 4$. It remains to determine $\pi_2(R)$ and
$\pi_3(R)$.

If $(E,h)$ is a pair over $S^3$,
then we have  $E=S^1\times S^3$ and $h=n(E,h)1_{S^1}\times o_{S^3}$ for a
well-defined integer  $n(E,h)\in
\Z$,
where $o_{S^3}\in H^3(S^3,\Z)$ is the canonical generator.
The bijection $P(S^3)\cong \Z$ given by $(E,h)\mapsto n(E,h)$
induces the isomorphism $\pi_3(R)\rightarrow \Z$ in view of
\ref{cla}.

Let us now consider a pair $(E,h)$ over $S^2$.
Note that $E$ is canonically oriented, in particular $H^3(E,\Z)=[E]\cdot\Z$.
Let $c\in H^2(S^2,\Z)$ be its first Chern class.
Then we define the tuple of integers
$$(k(E,h),n(E,h))=(<c,[S^2]>,<h,[E]>)\in \Z\oplus \Z\ .$$
The bijection
$P(S^2)\cong \Z\oplus \Z$  given by
$(E,h)\mapsto (k(E,h),n(E,h))$ defines the isomorphism
$\pi_2(R)\cong \Z\oplus \Z$ in view of \ref{cla}.\footnote{We leave it
  to the interested reader to check that these bijections are in fact
  homomorphisms.} \hB

\subsubsection{}
The computation of the homotopy groups of $R$ implies by the Hurewicz
theorem that
$H_0(R,\Z)=\Z$, $H_1(R,\Z)\cong 0$ and $H_2(R,\Z)\cong \Z\oplus \Z$.
By the universal coefficient theorem
$H^2(R,\Z)\cong \Z\oplus \Z$.
Recall that $\bc\in H^2(R,\Z)$ is the class of the projection
$\bc:R\rightarrow K(\Z,2)$. Let $\pi:\bE\rightarrow R$ be the
universal bundle and $\bh\in H^3(\bE,\Z)$ be the universal class.
\begin{ddd}
We define $\hat \bc:=-\pi_!(\bh)\in H^2(R,\Z)$.
\end{ddd}

\begin{lem}
We have $H^2(R,\Z)=\bc\Z\oplus \hat \bc\Z$.
\end{lem}
\proof
Using  the canonical isomorphisms $H^2(R,\Z)\cong \Hom(H_2(R,\Z),\Z)\cong
\Hom(\pi_2(R),\Z)$, where $x\in H^2(R,\Z)$ and $[f\colon S^2\to R]$ is
mapped to $<f^*x,[S^2]>$.
The identification $\pi_2(R)\cong \Z\oplus \Z$ above gives
$H^2(R,\Z)\cong \Z\oplus \Z$.
An inspection shows that
this isomorphism maps $a\bc+b\hat \bc$ to $(a,-b)$.
Therefore, $H^2(R,\Z)$ is freely generated by $\bc$ and $\hat \bc$.
\hB

\subsubsection{}
Let $\hat \bc$ be classified by a map
$\hat \bc:R\rightarrow K(\Z,2)$.
We will now determine the homotopy fiber $F$ of the map
$$(\bc, \hat \bc):R\rightarrow K(\Z,2)\times K(\Z,2)\ .$$

\begin{lem}
The homotopy fiber of $(\bc, \hat \bc)$ is
$K(\Z,3)$.
\end{lem}
\proof
We consider the long exact sequence of homotopy groups
$$\dots\rightarrow \pi_i(F)\rightarrow \pi_i(R)\rightarrow
\pi_i(K(\Z,2)\times K(\Z,2))\rightarrow\pi_{i-1}(F)\rightarrow\dots\ .$$
We immediately conclude that
$\pi_i(F)=0$ if $i\not\in \{1,2,3\}$.
Furthermore we see that
$\pi_3(F)\cong \pi_3(R)\cong \Z$.

Therefore the relevant part is now
$$0\rightarrow \pi_2(F)\rightarrow \pi_2(R){\xrightarrow{\pi_1(\bc,\hat\bc)}}
\pi_2(K(\Z,2)\times K(\Z,2))\rightarrow  \pi_1(F)\rightarrow 0\ .$$

Now we observe that
$(\bc,\hat \bc)$ induces an isomorphism in integral  cohomology of degree $i\le 2$.
Therefore it induces an isomorphism $\alpha:\pi_2(R)\cong
\pi_2(K(\Z,2)\times K(\Z,2))$. It follows that $\pi_i(F)=0$ for
$i\in\{1,2\}$.
\hB

\subsubsection{}
We now have seen that
$R$ is the total space of a bundle
$$K(\Z,3)\rightarrow R\stackrel{(\bc,\hat \bc)}{\rightarrow}
K(\Z,2)\times K(\Z,2)\ .$$
It remains to determine the invariant
$\kappa\in H^4(K(\Z,2)\times K(\Z,2),\Z)$ which determines this bundle.
To do this we compute the cohomology of $R$ up to degree $4$ and
then we determine the differential in the Serre spectral sequence of
the bundle.
We already know that
$$\begin{array}{|c|c|}\hline
n&H^n(R,\Z)\\\hline
0&\Z\\\hline
1&0\\\hline
2&\bc\Z\oplus \hat \bc\Z\\\hline\end{array}\ .$$

\subsubsection{}\label{sec:LKZ3}

We start with recalling the low-dimensional integral cohomology of
$LK(\Z,3)$. Note that $K(\Z,3)$ has the structure of an $H$-space
(because one possible model is $\Omega K(\Z,4)$), so
that $LK(\Z,3)$ is homotopy equivalent to $K(\Z,3)\times \Omega
K(\Z,3)$. Further note that $\Omega
K(\Z,3)\simeq K(\Z,2)$.
 We use that
$$\begin{array}{|c|c|c|}\hline n&H(K(\Z,2),\Z)&H(K(\Z,3),\Z)\\\hline
0&\Z&\Z\\\hline
1&0&0\\\hline
2&\Z&0\\\hline
3&0&\Z\\\hline
4&\Z&0\\\hline
5&0&0\\\hline
\end{array}\ .
$$
We now conclude by the K{\"u}nneth formula that
$$\begin{array}{|c|c|}\hline n&H(LK(\Z,3),\Z)\\\hline
0&\Z\\\hline
1&0\\\hline
2&\Z\\\hline
3&\Z\\\hline
4&\Z\\\hline
5&\Z\\\hline
\end{array}\ .$$

\subsubsection{}\label{sec:compute_E}
We compute the cohomology $H^3(R,\Z)$ using the Gysin sequence of
\begin{equation}\label{eq2}U(1)\rightarrow U\times LK(\Z,3)\rightarrow R\ .\end{equation}
Observe that
$$\begin{array}{ccc}
U\times LK(\Z,3)&\rightarrow&U\\
\downarrow&&\downarrow\\
R&\stackrel{\bc}{\rightarrow}&K(\Z,2)\end{array}$$
is a pull-back of $U(1)$-principal bundles.
Therefore the first Chern class of the $U(1)$-principal bundle
$U\times LK(\Z,3)\rightarrow R$ (with the \emph{diagonal $U(1)$-action}) is $\bc\in H^2(R,\Z)$.
We further use the fact that $U$ is contractible.
The relevant part of the Gysin sequence is
$$0\rightarrow H^3(R,\Z)\rightarrow H^3(LK(\Z,3),\Z)\rightarrow
H^2(R,\Z)\stackrel{\bc}{\rightarrow} H^4(R,\Z)\rightarrow
H^4(LK(\Z,3),\Z)\rightarrow H^3(R,\Z) \ .$$
Since $\bc$ is the first Chern class of $\pi\colon E\to R$, the above
principal bundle is isomorphic to $\bE$ and we can use the Gysin sequence for $\pi:\bE\rightarrow R$
$$\rightarrow H^3(\bE,\Z)\stackrel{\pi_!}{\rightarrow}
H^2(R,\Z)\stackrel{\bc}{\rightarrow} H^4(R,\Z)\rightarrow $$
to conclude that $\bc\cup \hat \bc=-\bc\cup \pi_!(\bh)=0$.
Therefore $\bc:H^2(R,\Z)\rightarrow H^4(R,\Z)$ is not injective.
Since $H^3(LK(\Z,3),\Z)\cong
\Z$ and $H^2(R,\Z)$ is free abelian this  implies that
$$H^3(R,\Z)=0\ .$$

\subsubsection{}
The map $\bc:R\rightarrow K(\Z,2)$ admits a natural split
$K(\Z,2)\rightarrow R$. It maps $x\in K(\Z,2)$ to
the class $ [u,\gamma]$, where $\gamma$ is the constant loop and
$u\in U_x$ is any point.
The split classifies the pair $(U,0)$ over $K(\Z,2)$.
The existence of the split implies that
$\bc$ generates a polynomial ring $\Z[\bc]$ as direct summand inside $H^*(R,\Z)$.

\subsubsection{}
In particular, $\bc^2\not=0$. Therefore the kernel
of $\bc:H^2(R,\Z)\rightarrow H^4(R,\Z)$ is generated by
$\hat \bc$. The Gysin sequence for (\ref{eq2}) now gives
$$0\rightarrow \Z\stackrel{\bc^2}{\rightarrow} H^4(R,\Z) \rightarrow
\Z\rightarrow 0\ ,$$
where the last copy $\Z$ is $H^4(LK(\Z,3),\Z)$.
This implies that $$H^4(R,\Z)\cong \bc^2\Z\oplus \Z\ .$$

We now show that $\bc^2$ and $\hat \bc^2$ generate $H^4(R,\Z)$ as a
$\Z$-module.
We consider the pair over $K(\Z,2)$ consisting of the trivial bundle
$\Pi:U(1)\times K(\Z,2)\rightarrow K(\Z,2)$
and the class $h=o_{U(1)}\times z\in H^3(U(1)\times K(\Z,2),\Z)$,
where $z\in H^2(K(\Z,2),\Z)$ is a generator.
This pair is classified by a map
$f:K(\Z,2)\rightarrow R$.
Let $F:U(1)\times K(\Z,2)\rightarrow \bE$ be defined by the pull-back
diagram
$$\begin{array}{ccc}
U(1)\times K(\Z,2)&\stackrel{F}{\rightarrow}&\bE\\
\Pi\downarrow&&\pi\downarrow\\
K(\Z,2)&\stackrel{f}{\rightarrow}&R
\end{array}\ .$$
Then we have
$f^*\bc=0$ and
$$f^*\hat \bc=-f^*\pi_!(\bh)=-\Pi_! F^*(\bh)=\Pi_!(h)=-z\ .$$
This shows that $\hat \bc\in H^2(R,\Z)$ generates a polynomial
ring isomorphic to $\Z[\hat \bc]$ inside $H^*(R;\Z)$. Furthermore, we see that $f^*(\hat
\bc^2)=z^2$ is primitive so that $\hat\bc^2$ must be primitive, too.
Thus $H^4(R,\Z)=\bc^2\Z\oplus \hat \bc^2\Z$.
Let us collect the results of our computations:
\begin{lem}
We have
$$\begin{array}{|c|c|}\hline n&H(R,\Z)\\\hline
0&\Z\\\hline
1&0\\\hline
2&\bc\Z\oplus \hat \bc\Z\\\hline
3&0\\\hline
4&\bc^2\Z\oplus \hat \bc^2\Z\\\hline
\end{array}\ .$$
\end{lem}

\subsubsection{}
We now finish the proof of Theorem \ref{main1}.
We consider the $E_4$-page of the Serre spectral sequence of the
fibration (\ref{eq3}).
$$\begin{array}{|c|c|c|c|c|c|}\hline
4&0&0&0&0&0\\\hline
3&\Z&0&*&0&*\\\hline
2&0&0&0&\Z&0\\\hline
1&0&0&0&0&0\\\hline
0&\Z&0&c\Z\oplus \hat c \Z&0&c^2\Z\oplus (c\cup\hat c)\Z\oplus  \hat c \Z\\\hline
*&0&1&2&3&4\\\hline
\end{array}\ .$$
We read off that
$$0\rightarrow \Z\stackrel{d^{0,3}_4}{\rightarrow} c^2\Z\oplus (c\cup\hat
c)\Z\oplus  \hat c \Z\rightarrow \bc^2\R\oplus \hat \bc^2 \Z\rightarrow
0\ .$$
The last map is the edge homomorphism and therefore induced by the map
$R\to K(\Z,2)\times K(\Z,2)$. Since under this map $c$ is mapped to
$\bc$ and $\hat c$ to $\hat \bc$,
$d^{0,3}_4$ is multiplication by $\pm c\cup \hat c$.
This finishes the proof of Theorem \ref{main1} \hB

\subsection{The $T$-transformation}

\subsubsection{}\label{red}
We have already observed that $\bE$ and $U\times LK(\Z,3)$ are
isomorphic $U(1)$-principal bundles over $R$, both having first Chern
class $\bc$. Since $U$ is contractible, we know the low dimensional cohomology of
$\bE$ by Section \ref{sec:LKZ3}. Using the Gysin sequence
of $\pi:\bE\rightarrow R$, we determine the generators in terms of
characteristic classes of $\bE$.
From
$$0\rightarrow H^0(R,\Z)\stackrel{\bc}{\rightarrow} H^2(R,\Z)\rightarrow
H^2(\bE,\Z)\rightarrow 0\ $$
we conclude that
$H^2(\bE,\Z)\cong \pi^*\hat \bc\Z$.
Finally, we get
$$0\rightarrow H^3(\bE,\Z)\stackrel{\pi_!}{\rightarrow}
H^2(R,\Z)\stackrel{\bc}{\rightarrow} H^4(R,\Z)\rightarrow
H^4(\bE,\Z)\rightarrow 0\ .$$
This shows that $H^3(\bE,\Z)\cong \Z$ and $H^4(\bE,\Z)\cong \pi^*\hat
\bc^2\Z$.
Since $\hat \bc=-\pi_!(\bh)$ generates the kernel of
$\bc:H^2(R,\Z)\rightarrow H^4(R,\Z)$ we have
$H^3(\bE,\Z)\cong \bh\Z$.
\begin{lem}\label{lem:cohomology_of_E}
$$\begin{array}{|c|c|}\hline n&H(\bE,\Z)\\\hline
0&\Z\\\hline
1&0\\\hline
2&\pi^* \hat \bc\Z\\\hline
3&\bh\Z \\\hline
4&\pi^*\hat \bc^2 \Z\\\hline
\end{array}\ .$$
\end{lem}

\subsubsection{}\label{rew}
The class $\hat \bc$ classifies a $U(1)$-principal bundle
$\hat \pi:\hat \bE\rightarrow R$. Since $\bc\cup\hat \bc=0$ and
$H^3(R,\Z)=0$ there exist unique classes
$h\in H^3(\bE,\Z)$ and $\hat h\in H^3(\hat \bE,\Z)$ such that
$(\bE,h)$ and $(\hat \bE,\hat h)$ are dual to each other, where we use
Corollary \ref{var}.

\begin{lem}
We have $h=\bh$.
\end{lem}
\proof
Let $r:\bV\rightarrow R$ denote the two-dimensional complex vector
bundle given by $\bV:=\bL\oplus \hat \bL$,
where $\bL$ and $\hat \bL$ are the hermitian line bundles associated
to $\bE$ and $\hat \bE$.
Then we can factor the associated unit sphere bundle as
$$S(\bV)\stackrel{s}{\rightarrow} P(\bV)\stackrel{t}{\rightarrow} R\ ,$$
where $P(\bV)$ is the projective bundle of $\bV$.
Let $\tilde c\in H^2(P(\bV),\Z)$ be the first Chern class of the
$U(1)$-principal bundle
$s:S(\bV)\rightarrow P(\bV)$.
By the Leray-Hirsch theorem $H^*(P(\bV),\Z)$ is a free module over
$H^*(R,\Z)$ generated by $1_{P(\bV)}\in H^0(P(\bV),\Z)$ and
$\tilde c$.
The line bundles $\bL$ and $\hat \bL$ induce two sections
$\bl,\hat \bl:R\rightarrow P(\bV)$ such that we have the following
pull-back diagrams
$$\begin{array}{ccccc}
\bE&\stackrel{\bi}{\rightarrow}&S(\bV)&\stackrel{\hat\bi}{\leftarrow}&\hat \bE\\
\pi\downarrow&&s\downarrow&&\hat \pi\downarrow\\
R&\stackrel{\bl}{\rightarrow}&P(\bV)&\stackrel{\hat \bl}{\leftarrow}&R
\end{array}\ .$$
Note that $\bl^*\tilde c=\bc$ and $\hat \bl^*\tilde c=\hat \bc$.
Since $t_!(\tilde c)=1=t_!\circ s_!(\Th)$
we have $s_!(\Th)=\tilde c+t^*b$ for some $b\in H^2(R,\Z)$.
This implies that
$\pi_!(h)=\pi_!\circ \bi^*(\Th)=l^*\circ s_!(\Th)=\bc+b$.
Analogously, we get
$\hat \pi_!(\hat h)=\hat \bc+b$.
Furthermore, we deduce from the projection formula that
$$\bc\cup(\bc+b)=\bc \cup\pi_!(h)=\pi_!(\pi^*\bc \cup h)=0,\quad  \hat \bc\cup (\hat \bc+b)=0\ .$$
Using the information about the ring structure of $H^*(R,\Z)$ it follows that
$\bc+b=m\hat \bc$ and $\hat \bc+b=n \bc$ for some $m,n\in \Z$.
Since $H^2(R,\Z)$ is freely generated by $\bc$ and $\hat \bc$ we
conclude that $m=n=-1$, i.e. $b=-(\hat \bc+\bc)$ so that
$\pi_!(h)=-\hat \bc$. By Lemma \ref{lem:cohomology_of_E}
each class in $x\in H^3(E,\Z)$ is a multiple of $\bh$.
Since $\pi_!(h)=\pi_!(\bh)$ we see that $h=\bh$.
\hB

\subsubsection{}

We also see that $\pi_!(\hat h)=-\bc$.
This shows that $\hat \bh:=\hat h\in H^3(\hat E,\Z)$ is a generator.
\begin{ddd}\label{dpa}
We define the dual universal pair to be $(\hat \bE,\hat \bh)$.
\end{ddd}
As in \ref{red} we have
\begin{kor}
$$\begin{array}{|c|c|}\hline n&H(\hat \bE,\Z)\\\hline
0&\Z\\\hline
1&0\\\hline
2&\pi^* \bc\Z\\\hline
3&\hat \bh\Z \\\hline
4&\pi^*\bc^2 \Z\\\hline
\end{array}\ .$$
\end{kor}

\subsubsection{}
Let $T:R\rightarrow R$ be the classifying map of the dual pair
$(\hat \bE,\hat \bh)$, covered by the $U(1)$-bundle map $T_E\colon
\hat\bE\to \bE$.

\begin{lem}\label{per4}
$T\circ T$ classifies $( E,\bh)$.
In particular, $T^2\sim  \id_R$.
\end{lem}
\proof
We have $T^*\bc=\hat \bc$. Furthermore,
$T^*\hat \bc=T^*\pi_!(-\bh)=\hat \pi_!(-T_E^*\bh)=\hat\pi_!(-\hat \bh)=\bc$.
Thus $(T\circ T)^*\bc=\bc$ and $(T\circ T)^*\hat \bc=\hat \bc$.
The underlying bundle of the pair classified by $T^2$ is
$ \pi: E\rightarrow R$.
Since $\pi_!(\bh)=-\hat \bc=-(T\circ T)^*\bc$
we must have $(T\circ T)^*(\bE,\bh)\cong (\bE,\bh)$.
\hB

\subsubsection{}
Recall from \ref{tr1} that $P(B)$ denotes the set of isomorphism classes of pairs
over $B$, and that we have a natural transformation of functors
$\Psi_B:[B,R]\rightarrow P(B)$.
The map $T:R\rightarrow R$ induces an involution
$T_*:[\dots,R]\rightarrow [\dots,R]$.

\begin{ddd}
We define the natural transformation of set-valued functors
$$T_{\dots}:P(\dots)\rightarrow P(\dots)$$ by
$$T_B:=\Psi_B\circ T_*\circ \Psi_B^{-1}\ .$$
We call it the \emph{$T$-duality transformation}.
\end{ddd}

\subsubsection{}

The following is a consequence of \ref{per4}.
\begin{kor}
Note that $T^2_B=\id$. In particular, the $T$-duality transformation
is an isomorphism of functors.
\end{kor}

\subsubsection{}

Let $(E,h)$ be a pair over $B$ and $c\in H^2(B,\Z)$ be the first Chern
class of $E$.
\begin{lem}
Any pair $(E,h)$ admits the  dual pair $(\hat E,\hat h)$ representing the
class $T_B([E,h])$. The first Chern class $\hat c\in H^2(B,\Z)$ of
$\hat E$
is given by $\hat c=-\pi_!(h)$. Furthermore, $c=-\hat \pi_!(\hat h)$.
\end{lem}
\proof
Let $f:B\rightarrow R$
classify the pair $(E,h)$. Then we let $(\hat E,\hat h)=f^*(\hat
\bE,\hat \bh)$. The relations between the Chern classes and the
three-dimensional cohomology classes follow from the corresponding
relations over $R$ obtained in \ref{rew}.
We have compatible pull-back diagrams
$$\begin{array}{ccc}
E&\rightarrow &\bE\\
\downarrow&&\downarrow\\
B&\rightarrow&R\end{array},\quad
\begin{array}{ccc}
\hat E&\rightarrow &\hat \bE\\
\downarrow&&\downarrow\\
B&\rightarrow&R\end{array}\ ,\quad
\begin{array}{ccc}
S(V)&\rightarrow &S(\bV)\\
\downarrow&&\downarrow\\
B&\rightarrow&R\end{array}\ .$$
We obtain the Thom class of $S(V)$ as a pull-back of the universal Thom class
of $S(\bV)$. Its restriction to $E$ and $\hat E$ gives $h$ and $\hat h$.
This shows that
$(E,h)$ and $(\hat E,\hat h)$ are in duality. \hB

\subsubsection{}
We consider pairs $(E,h)$ and $(\hat E,\hat h)$ over a space $B$.
Let $c,\hat c$ denote the first Chern classes of $E$ and $\hat E$.
\begin{lem}\label{lem2}
If $(E,h)$ and $(\hat E,\hat h)$ are dual to each other,
then we have $c=-\hat \pi_!(\hat h)$ and $\hat c=-\pi_!(h)$.
\end{lem}
\proof
Denote the canonical generators of the polynomial ring
$H^*(K(\Z,2)\times K(\Z,2),\Z)$ by $z,\hat z$ (instead of $c,\hat c$ -
we
do this in order to avoid notational
conflicts).Recall that we have a bundle
$$K(\Z,3)\rightarrow R\stackrel{(\bc,\hat \bc)}{\rightarrow}
K(\Z,2)\times K(\Z,2)$$
which is classified by $z\cup \hat z\in H^4(K(\Z,2)\times
K(\Z,2),\Z)$.
If $\bar f:B\rightarrow K(\Z,2)\times K(\Z,2)$ satisfies
$\bar f^*z\cup \bar f\hat z=0$, then it admits a lift
$f:B\rightarrow R$. Let
$\bar f:B\rightarrow K(\Z,2)\times K(\Z,2)$ be the classifying map
of the pair $(c,\hat c)$, i.e $\bar f^*z=c$ and $\bar f^*\hat z=\hat c$.
Then we have a lift $f:B\rightarrow R$.
Pulling back the universal pairs over $R$ we get
pairs $(E,h^\prime)$ and $(\hat E,\hat h^\prime)$ which are dual to
each other. Furthermore, $\pi_!(h^\prime)=-\hat c$ and $\hat
\pi_!(\hat h^\prime)=-c$. By \ref{var} we have
$h=h^\prime+\pi^*b$ and $\hat h=\hat h^\prime+\hat \pi^* b$ for some
$b\in H^3(B)$. Hence $\pi_!(h)=\pi_!(h^\prime)=-\hat c$ and $\pi_!(\hat
h)=\pi_!(\hat h^\prime)=-c$.\hB

\subsubsection{}
Note that there is a natural action of $H^3(B,\Z)$ on the set $P(B)$ given by
$\beta[E,h]:=(E,h+\pi^*\beta)$, $\beta\in H^3(B)$.
\begin{lem}
The $T$-duality transformation is equivariant with respect to this
action of $H^3(B,\Z)$.
\end{lem}
\proof
This is an immediate consequence of \ref{var}. \hB

\subsubsection{}

By Theorem \ref{uni} we already knew that the equivalence class of
pairs dual to $(E,h)$ is unique, if such  dual pairs exist at all.
The new information obtained from the study of the topology of the
classifying space is the existence of pairs dual to $(E,h)$.
More significantly, note that our proof of the uniqueness part of \ref{uni} involves Lemma
\ref{lem2}, whose proof also depends on the knowledge of the
topology of $R$.

\section{$T$-duality in twisted cohomology theories}

\subsection{Axioms of twisted cohomology}\label{fff}

\subsubsection{}
There may be many explicit models of a twisted cohomology theory which
lead to equivalent results, and examples abound in the literature. In
particular, this applies to the nature of a twist. What we will
describe here is a picture which should be the common core of the
various concrete realizations.

In any case the twists come as a pre-sheaf of pointed groupoids $B\mapsto
T(B)$ on the category of spaces. Let us fix some notation for the main
ingredients, which also recalls the concept of a pre-sheaf we use.
First of all $T(B)$ is a groupoid with a distinguished trivial object
$\theta_B$, giving rise to the trivial twist (i.e.~to no twist at
all). If $f:A\rightarrow B$ is a map of spaces, then there is a functor
$f^*:T(B) \rightarrow T(A)$ preserving the trivial twists.
Furthermore, if $g:B\rightarrow C$ is a second map,
then there exists a natural transformation
$$\Psi_{f,g}:f^*\circ g^*\rightarrow  (g\circ f)^*\ .$$
If $h:C\rightarrow D$ is a third map, then we require that
$$ \Psi_{f,h\circ g} \circ f^*\Psi_{g,h}=\Psi_{g\circ f,h}\circ h_*\Psi_{f,g}\
$$

\subsubsection{}\label{coup}
The following three requirements provide the coupling to topology.
\begin{enumerate}
\item
We require that there is a natural transformation
$c:T(\dots)\rightarrow H^3(\dots,\Z)$ (the latter is considered as a
pre-sheaf of categories
in a trivial manner, i.e.~with only identity morphisms) which classifies the isomorphism
classes of $T(B)$ for each $B$.
\item
If $\cH,\cH^\prime\in T(B)$ are equivalent objects, then we require that
$\Hom(\cH,\cH^\prime)$ is a
$H^2(B,\Z)$-torsor such that  the composition with
fixed morphisms gives isomorphisms of torsors.
Furthermore, we require that the torsor structure is
compatible with the pull-back. Note that we have natural bijections
$u:\Hom(\cH,\cH)\rightarrow H^2(B,\Z)$ which map compositions to sums.
\item
Let $\cK\in T(\Sigma(B\cup *))$,
where $\Sigma(B\cup *):=I\times B/(\{0\}\times B\cup \{1\}\times B)$
is the (reduced) suspension.
We have a homotopy
$h:I\times B\rightarrow \Sigma(B\cup *)$
from the constant map $p:B\mapsto *\into \Sigma(B\cup *)$ to itself given by
$h_t(b)=[t,b]$. It induces a morphism
$u(h):p^*\cK\rightarrow p^*\cK$ as explained in \ref{sec:hominv}.
We require that $u(u(h))$ and $c(\cK)$ correspond to each other under
the suspension isomorphism $H^3(\Sigma(B\cup *),\Z)\cong H^2(B,\Z)$.
\end{enumerate}

\subsubsection{}
 Let us list two examples.
\begin{enumerate}
\item
In our first example the objects of $T(B)$ are Hitchin gerbes.
Recall that a Hitchin gerbe over $X$ is a $U(1)$-extension
$\cH\rightarrow \cG$, where $\cG$ is an {\'e}tale  groupoid which represents the space $B$.
 A morphism in $T(B)$
is an equivalence class of equivalences of Hitchin gerbes $u:\cH\rightarrow
\cH^\prime$.
The isomorphism classes of Hitchin gerbes are
classified by the characteristic class $c(\cH)\in H^3(B,\Z)$.
We refer to \cite{bunkeschick04} for further details, in  particular
the torsor structure on the sets of morphisms.
\item
In the second example the objects of $T(B)$ are given by the set of continuous
maps $f:B\rightarrow K(\Z,3)$.
A morphism
$u:f\rightarrow f^\prime$
is then a homotopy class of homotopies from $f$ to $f^\prime$.
We set  $c(f):=[f]\in [B,K(\Z,3)]\cong H^3(B,\Z)$.
Recall that $LK(\Z,3)\simeq K(\Z,2)\times K(\Z,3)$. Let
$u:LK(\Z,3)\rightarrow K(\Z,2)$ be the first projection. The second
projection is given by the evaluation map $\ev_0$.
An automorphism of $f$ is a homotopy class $[\gamma]\in [B, LK(\Z,3)]$
with $\ev_0\circ \gamma=f$. Therefore, automorphisms are classified
by $[u\circ \gamma]\in [B,K(\Z,2)]\cong H^2(B,\Z)$.
\end{enumerate}

\subsubsection{}
In the following we fix some framework of twists and formulate
the axioms of a twisted cohomology theory in this framework.
We fix a  cohomology theory $h$ for which we want to define a twisted extension.

\begin{ddd}
A twisted cohomology theory $h$ extending $h$ associates to each
space $X$ and each twist  $\cH\in T(X)$ a $\Z$-graded
group $h(X,\cH)$. To
a map $f:Y\rightarrow X$ it associates a homomorphism
$$f^*:h(X,\cH)\rightarrow h(Y,f^*\cH)\ .$$
To a morphism $u:\cH\rightarrow
\cH^\prime$ of twists it associates
an isomorphism, natural with respect to pullbacks,
 $$u^*:h(X,\cH)\rightarrow h(X,\cH^\prime)\ .$$
Finally, we require an integration map
$$p_!:h(Y,p^*\cH)\rightarrow h(X, \cH)$$
of degree $\dim(Y)-\dim(X)$
for a proper  $h$-oriented map $p:Y\rightarrow X$. Integration shall
be natural with respect to morphisms in $T(X)$.
\end{ddd}
These structures must satisfy the axioms described below.

\begin{axiom}[Extension]
Let $\theta_X\in T(X)$ denote the trivial twist.
There exists a canonical isomorphism
$h(X,\theta_X)\rightarrow h(X)$
which preserves pull-back and integration over the
fiber.
\end{axiom}

\begin{axiom}[Functoriality]\label{func}
If $g:Z\rightarrow Y$ is a second map,
then we have
$$\Psi_{g,f}(\cH)^* \circ (f\circ g)^*=g^*\circ f^*\ .$$
If $v:\cH^{\prime\prime}\rightarrow \cH$ is another morphism of
twists, then we have
$$v^*\circ u^*=(u\circ v)^*\ .$$
\end{axiom}

\subsubsection{}
\label{sec:hominv}

Assume that $h:\R\times
Y\rightarrow  X$ is a homotopy from $f_0$ to $f_1$,
i.e. $f_k=i_k^*(h)$,
where $i_k:Y\rightarrow \R\times Y$ is given by $i_k(x)=(k,x)$,
$k=0,1$. Define $F\colon \R\times Y\to \R\times X; (t,y)\mapsto
(t,h(t,y))$. Observe that for $\cH\in T(X)$ the twists $(\id_\R\times f_0)^*\pr_2^*\cH$
and
$F^*\pr_2^*\cH$ on $\R\times Y$ are isomorphic, since we can by
assumption read off the isomorphism class from the pullbacks of the
corresponding classifying cohomology class, which are equal by
homotopy invariance of cohomology. We define
\begin{equation*}
u(h):(\id_\R\times f_0)^*\pr_2^*\cH\rightarrow F^*\pr_2^*\cH
\end{equation*}
to be the unique
morphism of twists such that
$$f_0^*\cH\stackrel{can}{\cong}i_0^*\circ (\id_\R\times f_0)^*\circ \pr_2^*\cH
\stackrel{i_0^*(u(h))}{\cong}
i_0^*\circ F^*\circ \pr_2^*\cH\stackrel{can}{\cong} f_0^*\cH$$
is the identity. The morphism $u(h)$ is determined uniquely this way since
$i_0^*\colon H^2(\reals\times Y,\Z)\to H^2(Y,\integers)$
is an isomorphism. The canonical isomorphisms are induced by Axiom
\ref{func}. Note that $u(h)$ is natural with respect to morphisms in
$T(X)$.

Finally we define
$$v(F):f_0^*\cH\stackrel{can}{\cong} i_1^*\circ (\id_\R\times f_0)^*\circ
\pr_2^*\stackrel{i_1^*(u(h))}{\cong} i_1^*\circ F^*\circ \pr_2^*\cH\stackrel{can}{\cong}
 f_1^*\cH$$
\begin{axiom}[Homotopy invariance]\label{homaxiom}
With these conventions we require that $$v^*\circ f_1^*=f_0^*\ .$$
\end{axiom}


\begin{axiom}[Integration]
\mbox{}\begin{enumerate}
\item{\bf Functoriality}
If $q:Z\rightarrow Y$ is a further proper $h$-oriented map, then we
have
$$p_!\circ q_!\circ  \Psi_{p,q}(\cH)^*
=(q\circ p)_!:h(Z,(q\circ p)^*(\cH))\rightarrow
h(X,\cH)\ .$$
\item{\bf Naturality}
If $g:Z\rightarrow X$ is a further map, then we
have the Cartesian diagram
$$\begin{array}{rcr}
Z\times_XY&\stackrel{p^*g}{\rightarrow}&Y\\
g^*p\downarrow&&p\downarrow\\
Z&\stackrel{g}{\rightarrow}&X
\end{array}\ ,$$ and we require that
$$(g^*p)_!\circ (\Psi_{g,g^*p})(\cH)^*\circ (\Psi_{p,p^*g}(\cH)^*)^{-1}\circ (p^*g)^*
=g^*\circ p_!\ .
$$
\end{enumerate}
\end{axiom}

\begin{axiom}[Mayer-Vietoris sequence]
 If $X=U\cup V$ is a decomposition   by open subsets, then we can find a function
$\phi:X\rightarrow \R$ such that $\phi_{|X\setminus U}=1$,
$\phi_{X\setminus V}=-1$, and the inclusion
$i\colon (Y:=\{\phi=0\})\rightarrow X$ is a proper naturally $h$-oriented map.
Let $j:Y\rightarrow U\cap V$,
$g:U\rightarrow X$, $h:V\rightarrow X$, $k:U\cap V\rightarrow U$,
$l:U\cap V\rightarrow V$, and $r:U\cap V\rightarrow X$  denote the
inclusions, and define $\delta:=i_!\circ j^*$\footnote{Note that $\delta$ is
independent of the choice of $\phi$.}.
Then we require that the following sequence is exact:
$$\dots \rightarrow h(U\cap V,r^*\cH)\stackrel{\delta}{\rightarrow}
h(X,\cH)\stackrel{(g^*,h^*)}{\rightarrow} h(U,g^*\cH)\oplus
h(V,h^*\cH)\stackrel{k^*-l^*}{\rightarrow} h(U\cap
V,r^*\cH)\rightarrow \dots \ ,$$
where some canonical isomorphisms are suppressed in the notation.
\end{axiom}

\subsubsection{}
Examples of twisted cohomology theories which satisfy these axioms
(on the category of smooth manifolds and smooth maps)
are twisted de Rham cohomology and twisted $Spin^c$-cobordism theory
\cite{bunkeschick04} and \cite{bunkeschick041}.
In these examples twists are Hitchin gerbes. As indicated in
\cite{bunkeschick041} there should also be a twisted version of
complex
$K$-theory. In this case the missing piece in the literature is a nice
description of integration over the fiber and also of the boundary operator in the
Mayer-Vietoris sequence.

\subsection{$T$-admissibility}

\subsubsection{}\label{poi}
We consider the unit sphere $S\subset \C^2=\C\oplus \C$.
Let $E:=S^1$ and $\hat E:=S^1$. We consider the embeddings
$i:E\rightarrow S$, $i(z)=(z,0)$ and $\hat i:\hat E\rightarrow S$,
$\hat i(\hat z)=(0,\hat z)$. Let $T:=E\times \hat E$ and
$p:T\rightarrow E$ and $\hat p:T\rightarrow \hat E$ denote the
projections.
We define the homotopy  $h:I\times T\rightarrow S$ from $i\circ p$ to
$\hat i\circ \hat p$ by
$$h_t(z,\hat z):=\frac{1}{\sqrt{2}}(\sqrt{1-t^2}z,t\hat z)\ .
$$
Let $\cK\in T(S)$ be a twist such that $<c(\cK),[S]>=1$.
We define $\cH:=i^*\cK$ and $\hat \cH:=\hat i^*\cK$.
The homotopy $h$ induces a unique morphism $$u:\hat p^* \hat\cH=\hat
p^*\hat i^*\cK\stackrel{\Psi_{\hat p,\hat i}(\cK)}{\cong} (\hat i\circ \hat p)^*\cK\stackrel{u(h)}{\cong}
(i\circ p)^*\cK\stackrel{\Psi_{p,i}(\cK)^{-1}}{\cong} p^*i^*\cK=p^*\cH\ ,$$
where $u(h)$ is defined in Section \ref{sec:hominv}.

\subsubsection{}
Let $h$ be a twisted cohomology theory. Note that
$\hat p$ is canonically $h$-oriented since $T\hat E$ is canonically
trivialized by the $U(1)$-action.
\begin{ddd}
We say that the twisted cohomology theory $h$ is $T$-admissible if
$$\hat p_!\circ u(h)^*\circ p^*:h(E,\cH)\rightarrow h(\hat E,\hat \cH)$$
is an isomorphism. Note that the map has degree $-1$.
\end{ddd}

\subsubsection{}
Naturality implies that $T$-admissibility does not depend on the
choice of $\cK$ inside its isomorphism class.

\subsubsection{}
We show now how one can check $T$-admissibility in practice.

Let $S/(E\cup \hat E)$ be the quotient space of $S$ where $i(E)$ and
$\hat{i}(\hat E)$ are identified to one point. We have a natural
identification $r:S/(E\cup \hat E)\cong \Sigma (T\cup
*)$ given by the homotopy $h$ used in Section \ref{poi}. Note that $p^*r^*:H^3( \Sigma (T\cup
*),\Z)\rightarrow H^3(S,\Z)$ is an isomorphism, where $p\colon S\to
S/(E\cup \hat E)$ is the projection.
Thus, we  can choose $\cK:=p^*r^*\tilde \cK$ for some twist
$\tilde \cK\in T(\Sigma (T\cup
*))$. Note that $c(\cK)\in H^3(\Sigma (T\cup
*),\Z)\cong \Z$ is a generator.
Since $H^3(*,\Z)=0=H^2(*,\Z)$,  the restriction of $\tilde \cK$ to
the base point is the trivial twist. Then we obtain canonical
morphisms
$\cH\cong \theta_E$ and $\hat \cH \cong \theta_{\hat E}$.
The homotopy $r\circ h$ induces now a canonical morphism
$u(r\circ h):\theta_T\rightarrow \theta_T$. By the third property
stated in \ref{coup} we know that $u(u(r\circ h))\in H^2(T,\Z)\cong
\Z$ is a generator, too. The determination of this generator involves
the precise understanding of the isomorphism $p^*r^*$ and of the
suspension isomorphism.

Note that $H^2(T,\Z)$ acts naturally on $h(T)$ via the identifications
$H^2(T,\Z)\cong \Hom(\theta_T,\theta_T)$ and $h(T)\cong
h(T,\theta_T)$.
For $g\in H^2(T,\Z)$ we denote this action by $g^*$.

Therefore, in order to check that the cohomology theory $h$ is
$T$-admissible, it suffices to show
that
$$\hat p_!\circ g^* \circ p^*:h(E)\rightarrow h(\hat E)$$
is an isomorphism if $g\in H^2(T,\Z)$ is a generator.

\subsubsection{}

\begin{lem}
Twisted $K$-theory is $T$-admissible.
\end{lem}
\proof
Let $l\in K^0(T)$ be the class of the line bundle over $T$
with first Chern class equal to $g\in H^2(T,\Z)\cong \Z$. Then
$g^*$ is induced by the cup product with $l$.
Let $1\in K^0(S^1)$ and $u\in K^1(S^1)$ be the generators.
One can compute
\begin{eqnarray*}
\hat p_!\circ g^* \circ p^*(1)&=&g\ B(u)\\
\hat p_!\circ g^* \circ p^*(u)&=&1\ ,
\end{eqnarray*}
where $B:K^{1}\rightarrow K^{-1}$ is the Bott periodicity transformation.
This is indeed an isomorphism if $g\in\{1,-1\}$. \hB

\subsubsection{}
We consider the graded ring $R:=\R[z,z^{-1}]$ where $\deg(z)=-2$ and
the twisted cohomology $H_R(X,\cH)$, where we use $z$ in order to couple the twist.

\begin{lem}
Twisted cohomology with coefficients in $R$ is $T$-admissible.
\end{lem}
\proof
The action of $g\in H_R^2(T,\Z)$ is given by the cup-product with
$1+z g_\R$, where $g_\R$ is the image of $g$ in $H_\R(T)$.
By a simple computation
\begin{eqnarray*}
\hat p_!\circ g^* \circ p^*(1)&=&z g\ u\\
\hat p_!\circ g^* \circ p^*(u)&=&1
\end{eqnarray*}
This is indeed an isomorphism if $g\not=0$. \hB

\subsubsection{}
$T$-admissibility is a strong condition on $h$.
It implies for example that
$$p_!\circ g^*\circ \hat p^*\circ \hat p_!\circ g^*\circ p^* :h(E)\rightarrow h(E)$$
is an isomorphism of degree $-2$. This isomorphism induces a
two-periodicity of $h(E)$. Here is a non-example.
\begin{lem}
Twisted $Spin^c$-cobordism is not $T$-admissible.
\end{lem}
\proof
$\MSpin^c(S^1)$ is not two-periodic since it is concentrated in degree
$\le 1$.
\hB


\subsection{$T$-duality isomorphisms}

\subsubsection{}
We consider two pairs $(E,h)$ and $(\hat E,\hat h)$ over $B$ which are
dual to each other. We use the notation of \ref{dua}.
Let $\Th\in H^3(S(V),\Z)$ be a Thom class.
We choose a twist $\cK\in T(S(V))$ such that $c(\cK)=\Th$.
Then we define $\cH:=i^*\cK\in T(E)$ and $\hat \cH:=\hat i^*\cK\in
T(\hat E)$. We have $c(\cH)=h$ and $c(\hat \cH)=\hat h$.
We consider the diagram
\begin{equation}\label{cd}\begin{array}{ccccc}
&&E\times_B\hat E&&\\
&p\swarrow&&\hat p\searrow&\\
E&&q\downarrow&&\hat E\\
&\pi \searrow&&\hat \pi \swarrow\\
&&B&&
\end{array}
\ .\end{equation}
This is the parameterized version of the situation considered in
\ref{poi}. In particular, we have a homotopy
$h:I\times E\times_M\hat E\rightarrow S(V)$ from $i\circ p$ to $\hat i\circ \hat p$.
It induces the morphism
 $$u:\hat p^* \hat\cH=\hat
p^*\hat i^*\cK\stackrel{\Psi_{\hat p,\hat i}(\cK)}{\cong} (\hat i\circ \hat p)^*\cK\stackrel{u(h)}{\cong}
(i\circ p)^*\cK\stackrel{\Psi_{p,i}(\cK)^{-1}}{\cong}
p^*i^*\cK=p^*\cH$$
which is natural under pullback of bundles.

Let $h$ be a twisted cohomology theory.
\begin{ddd}\label{trd}
We define the $T$-duality transformation
$$T:=\hat p_!\circ u^*\circ p^*:h(E,\cH)\rightarrow h(\hat E,\hat
\cH)\ .$$
\end{ddd}

\subsubsection{}
The main theorem of the present section is the following.
Assume that $B$ is  homotopy equivalent to a finite complex.
\begin{theorem}\label{main2}
If $h$ is $T$-admissible, then the
$T$-duality transformation $T$ is an isomorphism.
\end{theorem}
\proof
Let $f:A\rightarrow B$ be a map. Then we use the pull-back of $\cK$ in
order to define the duality transformation $T$ over $A$.
Let $F:f^*E\rightarrow E$ and $\hat F:f^*\hat E\rightarrow \hat E$ be
the induced maps. The statement of the following Lemma involves
various (not explicitly written) canonical isomorphisms.
\begin{lem}\label{l2}
We have
$$T\circ F^*=\hat F^* \circ T:h(E,\cH)\rightarrow h(f^*\hat E,\hat
F^* \hat \cH)\ .$$
\end{lem}
Assume that we have a decomposition $B=U\cup V$ with open subsets $U$
and $V$ and let $j:U\cap
V\rightarrow B$ denote the inclusion.
By taking pre-images with respect to $\pi$ and $\hat \pi$
we obtain associated decompositions
$E=E_U\cup E_V$ and $\hat E=\hat E_U\cup \hat E_V$.
Let $f:E_U\cap E_V\rightarrow E$ and $\hat f:\hat E_U\cap \hat
E_V\rightarrow \hat E$ denote the inclusions. Finally let
$\delta:h(E_U\cap E_V,f^*\cH)\rightarrow
h(E,\cH)$ and
$\hat \delta :h(\hat E_U\cap \hat
E_V,\hat f^*\hat \cH )\rightarrow h(\hat E,\hat \cH)$ denote the boundary operators in the
Mayer-Vietoris sequences.
\begin{lem}\label{l3}
We have
$$T\circ \delta=\hat \delta\circ T:h(E_U\cap E_V,f^*\cH)\rightarrow
h(\hat E,\hat \cH)\ .$$
\end{lem}
Assuming these lemmas, the proof of the theorem now goes by induction
on the number of cells of $B$. The induction starts with any
contractible base since $h$ is
$T$-admissible, using naturality and homotopy invariance.  In the
induction step we adjoin a cell. We use Lemma \ref{l2} and \ref{l3}
in order to see that $T$ induces a map of Mayer-Vietoris sequences.
The induction step now follows from the five-lemma.
\hB

\subsubsection{}
We now prove Lemma \ref{l2}. Let $G:f^*E\times_Af^*\hat E\rightarrow
E\times_B\hat E$ be the induced map.
The assertion follows from the following computation, omitting a
number of canonical isomorphisms.
\begin{eqnarray*}
\hat F^*\circ T&=&\hat F^*\circ \hat p_!\circ u^*\circ p^*\\
&=&(f^*p)_!\circ G^*\circ u^*\circ p^*\\
&=&(f^*p)_!\circ (G^*u)^* \circ G^*\circ p^*\\
&=&   (F^*p)_!\circ (G^*u)^* \circ (F^*p)^*\circ F^*\\
&=&T\circ F^*\ . 
\end{eqnarray*}
\hB
\subsubsection{}
We now prove Lemma \ref{l3}.
Let $\phi\in C(B)$ be a function which takes the value
$-1$ on $B\setminus V$ and $1$ on $B\setminus U$, and such
that
the inclusion of $i\colon (Y:=\{\phi=0\})\rightarrow B$
is canonically $h$-oriented.
 Let  $k:Y\rightarrow U\cap V$,  $I:=\pi^*i$, $\hat I:=\hat \pi^*i$, $K:=\pi^*k$, and
$\hat K:=\hat \pi^*k$ denote the
corresponding inclusions.
Note that $I$ and $\hat I$ have a trivialized normal bundle and
thus are canonically $h$-oriented. We have
$\delta=I_!\circ K^*$ and $\hat \delta=\hat I_!\circ \hat K^*$.
Furthermore, we set $\tilde I:=q^*i$ and $\tilde K=q^*k$.
Finally let
$J:=\pi^*j,\hat J:=\hat \pi^*j$, and
$G:=q^*j$ denote the corresponding embeddings over $j\colon U\cap
V\into B$.
The assertion of the Lemma now follows from the following computation,
where canonical isomorphisms are omitted.
\begin{eqnarray*}
\hat \delta \circ T&=& \hat I_!\circ \hat K^*\circ (\hat J^*p)_!\circ   (G^*u)^* \circ
(J^*p)^*\\
&=& \hat I_! \circ (\hat I^*p)_!\circ (\tilde I^*u)^*\circ (I^*p)^*
\circ K^*\\
&=& \hat p_!\circ\tilde I_! \circ (\tilde I^*u)^*\circ (I^*p)^*
\circ K^*\\
&=& \hat p_!\circ u^* \circ\tilde I_!\circ (I^*p)^* \circ K^*\\
&=& \hat p_!\circ u^* \circ p^*\circ I_!\circ K^*\\
&=& T\circ \delta \ .
\end{eqnarray*}
\hB

\section{Examples}\label{ex}

\subsection{The computation of twisted $K$-theory for $3$-manifolds}\label{dre}

\subsubsection{}
If $E$ is a closed oriented $3$-manifold then isomorphism classes of twists $\cH$ on $E$ are
classified by the number $<c(\cH),[E]>\in \Z$.
We fix an equivalence class of twists corresponding to $n\in
\Z$. Representatives can be pulled back from $S^3$ using a map of degree $1$.
Note that $K(E,\cH)$ is independent of the twist in its class up to a non-canonical
isomorphism. In the present subsection we want to compute the
isomorphism class of this group which we will denote by $K(E,n)$.
Our computation is based on the Mayer-Vietoris sequence.

\subsubsection{}
We choose a ball  $U\subset E$. Then we have a decomposition $E=U\cup V$ such
that $U\cap V\sim S^2$. We identify the twists on $U$ and $V$ with
the trivial twist. We can arrange that under the degree $1$ map to
$S^3$ the set $U$ is mapped to the complement of the south pole and $V$ is
mapped to the complement of the north pole.

Using the relation between twists and morphisms in the suspension $S^3$
of $S^2$ and naturality, we see that a twist in the class $n$ is given by the transition
morphism $v:\theta_{S^2}\rightarrow \theta_{S^2}$
such that $<u(v),[S^2]>=\pm n$. Let
$u^*:K(S^2)\rightarrow K(S^2)$ denote the corresponding automorphism.
It acts by the cup product with the class of the line bundle of
degree $\pm n$.
Then the Mayer-Vietoris sequence reads
$$\rightarrow K(S^2)\stackrel{\delta}{\rightarrow} K(E,n)\rightarrow K(U)\oplus
K(V)\stackrel{u^*\circ i^*-j^*}{\rightarrow} K(S^2)\rightarrow \ ,$$
where $i:S^2\rightarrow U$ and $j:S^2\rightarrow V$ are the
inclusions. At this point we have fixed the sign of the class of the twist.

\subsubsection{}
We identify $K(S^2)\cong \Z I\oplus  \Z \theta$,
where $I$ is represented by the trivial one-dimensional bundle,
and $\theta$ is represented by the difference of a line bundle of
degree one and a trivial line bundle. We then have
$u^*I=I\pm n\theta$ and $u^*\theta=\theta$.

\subsubsection{}
The Mayer-Vietoris sequence gives
$$0\rightarrow K^0(E,n)\rightarrow \Z\oplus
K^0(V)\xrightarrow{a:=u^*\circ i^*-j^*}
\Z I\oplus \Z\theta \stackrel{\delta}{\rightarrow} K^1(E,n)\rightarrow
K^1(V)\rightarrow 0\ .$$
The restriction of $a$ to the first summand maps $1\in \Z$ to
$I\oplus \pm n\theta$. If $x\in K^0(V)$,
then $x_{|S^2}=\dim(x) I+ <c_1(x)_{|S^2},[S^2]>\theta$. Now observe that
$<c_1(x)_{|S^2},[S^2]>=0$ since $S^2$ bounds in $V$.
Therefore we have $a(k,x)=(k-\dim(x))I\pm kn\theta$.
We conclude that for $n\ne 0$
$$K^0(E,n)\cong \tilde K^0(V)\cong \tilde K^0(E)\ ,$$
where $\tilde K^0(E):=\ker(\dim)$ is the reduced group.
Furthermore, $K^1(E,n)$ fits into a sequence
$$0\rightarrow \Z/n\Z\rightarrow K^1(E,n)\rightarrow K^1(V)\rightarrow
0\ .$$
Note that $K^1(V)$ is free abelian and satisfies
$\rk\:K^1(V)= \rk\:K^1(E)-1$. In particular we get
$$K^1(E,n)\iso\Z^{\rk\:K^1(E)-1}\oplus\Z/n\Z\ .$$

\subsubsection{}\label{eki}
Let $M$ be a closed oriented surface of genus $g$.
The $U(1)$-principal bundles over $M$ are classified by
the first Chern class. Let $\pi:E_k\rightarrow M$ be the bundle
with first Chern class $<c(E_k),[M]>=k$.
We use the Gysin sequence in order to compute the integral cohomology
of $E_k$.
We get
$$\begin{array}{|c|c|c|}\hline
i&H(E_k,\Z), k\not=0&H(E_0,\Z)\\\hline
0&\Z&\Z\\\hline
1&\Z^{2g}&\Z^{2g+1}\\\hline
2&\Z^{2g}\oplus \Z/k\Z&\Z^{2g+1}\\\hline
3&\Z&\Z\\\hline
\end{array}\ .$$

\subsubsection{}
We now compute the $K$-theory of $E_k$ using the Atiyah-Hirzebruch
spectral sequence.
The second page in the case $k\not=0$ looks like (vertically periodic)
$$\begin{array}{|c|c|c|c|c|c|}\hline
2&\Z&\Z^{2g}&\Z^{2g}\oplus \Z/k\Z&\Z\\\hline
1&0&0&0&0\\\hline
0&\Z&\Z^{2g}&\Z^{2g}\oplus \Z/k\Z&\Z\\\hline
&0&1&2&3\\\hline
\end{array}\ .$$
The only possibly non-trivial differential is $d_3^{2,3}$.
But since the spectral sequence degenerates rationally, the differential is
trivial. We get
$$\begin{array}{|c|c|c|}\hline
i&K(E_k), k\not=0&K(E_0)\\\hline
0&\Z^{2g+1}\oplus \Z/k\Z&\Z^{2g+2}\\\hline
1&\Z^{2g+1}&\Z^{2g+2}\\\hline
\end{array}\ .$$

\subsubsection{}
We now use this result in order to compute $K(E_k,n)$. We get for $n\not=0$
$$\begin{array}{|c|c|c|}\hline
i&K(E_k,n),k\not=0&K(E_0,n)\\\hline
0&\Z^{2g}\oplus \Z/k\Z&\Z^{2g+1}\\\hline
1&\Z^{2g}\oplus \Z/n\Z&\Z^{2g+1}\oplus \Z/n\Z\\\hline
\end{array}\ .$$

\subsubsection{}
Let us now verify that this computation confirms $T$-duality.
In fact, the unique dual pair of $(E_k,no_{E_k})$ is
$(E_n,-ko_{E_n})$. Thus $T$-duality predicts an isomorphism
$$K(E_k,n)\cong K(E_{-n},-k)$$ of degree $-1$. This is in fact compatible
with the results above.

\subsection{Line bundles over $\C P^r$}\label{sec:line-bundles-over}

\subsubsection{}
Let $p_n:E_{n,r}\rightarrow \C P^r$ be the $U(1)$-principal bundle with first
Chern class $nz$, where $z\in H^2(\C P^r,\Z)$ is the canonical generator.
We first compute $H(E_n,\Z)$ using the Gysin sequence for $p_n$.
We get
$$\begin{array}{|c|c|c|}\hline k&H^k(E_{n,r},\Z),n\not=0&H^k(E_{0,r},\Z)\\\hline
0&\Z&\Z\\\hline
2l, 1\le l\le r&\Z/n\Z&\Z\\\hline
2l+1, 1\le l\le r-1&0&\Z\\\hline
2r+1&\Z&\Z\\\hline
\end{array}\ .$$

Note that $r=\infty$ is permitted in the construction and calculation,
and that $E_{n,\infty}$ is a model for $B\Z/n\Z$.

\subsubsection{}
We compute the $K$-theory of $E_{n,r}$ using the Atiyah-Hirzebruch
spectral sequence. 
 We observe that this sequence degenerates. We get
$$\begin{array}{|c|c|c|}\hline &K(E_{n,r}), n\not=0&K(E_{0,r})\\\hline
0&\Z\oplus A_{n^r}&\Z^{2r+1}\\\hline
1&\Z&\Z^{2r+1}\\\hline
\end{array}\ ,$$
where $A_{n^r}$  is an abelian group with $n^r$ elements and with
composition series with subquotients $\Z/n\Z$. Using Atiyah's
completion theorem $\lim\limits_{\leftarrow_r}
\tilde K^0(E_{n,r})\iso \widehat{R(G)}$ we get extra information about
these groups, e.g.~that the limit is torsion-free. In the particularly
simple case $n=2$ we have
$\widehat{R(G)}\iso \Z_{(2)}$, which implies that $A_{2^r}\iso
Z/2^r\Z$ is cyclic. For other $n$, in particular if $n$ is a prime
number, $A_{n^r}$ can also be computed explicitly by looking at the
completion theorem and suitable Leray-Serre spectral sequences; we
leave this as an exercise to the reader. A precise answer can be found
e.g. in the book of Gilkey \cite{gilkey}, Thm. 4.6.7.

\subsubsection{}
The computation of the cohomology shows that for $r>1$ only $E_0$ admits
non-trivial twists (the case $r=1$ is covered by 
Section \ref{dre}). Let us fix the generator $g\in H^3(E_0)$
such that $(p_{0})_{!}(g)=z$.
Then twists $\cH$ over $E_0$ are classified by an integer $k\in \Z$
such that $c(\cH)=kg$.
 Let $K(E_0,k)$ be the isomorphism class of the twisted
$K$-theory for the twists in the class $k\in \Z$.
We can now apply $T$-duality in order to compute this group.
In fact, the unique dual pair of
$(E_0,k)$ is $(E_{-k},0)$. Thus we get
$$\begin{array}{|c|c|}\hline &K(E_0,k)\\\hline
0&\Z\\\hline
1&\Z\oplus A_{n^r}\\\hline
\end{array}\ .$$

Note that the calculations of this section, using the results of the present paper, rely on the fact that twisted K-theory is
a twisted cohomology theory in the sense of our axioms. As explained
earlier, no complete account of such a theory seems to be available in
the literature
.
\subsection{An example where torsion plays a role}\label{sec:an-example-where}

\subsubsection{}
As the base space we consider the total space  of the bundle
$p_k:E_{k,r}\rightarrow M$ as in Section \ref{sec:line-bundles-over}
for a prime number $k>1$ and for $r>1$, i.e. we set $B:=E_{k,r}$.
We fix a class $0\ne c\in H^2(B,\Z)$ and let $F_c$ denote the corresponding
$U(1)$-principal bundle over $B$. Since $c$ generates
$H^*(B,\Z)$ as a ring, except for the top degree, the Gysin sequence
of $F_c$ shows that its cohomology vanishes in degrees $1<i<2r+1$, and
$H^1(F_c,\Z)\iso\Z$.

Choose therefore $0=h\in H^3(F_c,\Z)$. Since $H^3(B,\Z)=0$, there is a
unique dual pair $(F_{\hat{c}},\hat h)$ with Chern class $\hat{c}=-\pi_!(h)=0$ and
$\hat{h}\in H^3(F_{\hat{c}},\Z)$. Since $\hat{c}=0$, $F_{\hat c}$ is the
trivial bundle, therefore its cohomology is
$H^l(F_{\hat{c}},\Z)\cong H^l(B,\Z)\oplus H^{l-1}(B,\Z)$,
e.g.~$H^3(F_{\hat{c}},\Z)\cong H^2(B,\Z)\cong \Z/k\Z$. Now $\hat{h}$ is the unique
  class with $\hat\pi_!(\hat h)=-c$, i.e.~corresponds to $-c$ under the
  isomorphism $H^3(F_{\hat c},\Z)=H^2(B,\Z)$. Clearly, if we only
  worked with differential forms as is done in \cite{bem}, then we
  could not distinguish this torsion twist from the trivial one.


\subsubsection{}
The 
 Atiyah-Hirzebruch spectral sequence for $K(F_c)$ degenerates. This  shows that
  $K_0(F_c)\iso \Z^2\iso K_1(F_c)$, whereas $K_0(F_0)\iso
  K_0(B)\oplus K_1(B)\iso \Z\oplus \Z \oplus A_{k^r}\iso
  K_1(F_0)$.  The $T$-duality
  isomorphism identifies $K(F_c)$ with $K(F_0,\hat h)$ for the torsion
  twist $\hat h$. In particular we see that $K(F_0)\not\cong K(F_0,\hat
  h)$ which shows that the torsion part of the twist is important.

\subsection{Iterated $T$-duality}
\label{sec:iterated-t-duality}

\subsubsection{}
Let $T$ denote the group $U(1)\times U(1)$.
\begin{ddd}
Two principal $T$-bundles $F\rightarrow B$ and $F^\prime\rightarrow B$
are isomorphic if there exists an isomorphism of fiber bundles
$$\begin{array}{ccc}F&\stackrel{U}{\rightarrow}&F^\prime\\
\downarrow&&\downarrow\\
B&=&B\end{array} $$
such that $U$ is $T$-equivariant.\end{ddd}

\subsubsection{}
The group of automorphisms of $T$ is $GL(2,\Z)$.
If we identify $T\cong \R^2/\Z^2$, then the action
of this group on $T$ is induced by the linear action on $\R^2$.
Let $\phi\in GL(2,\Z)$.
\begin{ddd}
Two principal $T$-bundles $F\rightarrow B$ and $F^\prime\rightarrow B$
are $\phi$-twisted isomorphic if there exists an isomorphism of fiber bundles
$$\begin{array}{ccc}F&\stackrel{U}{\rightarrow}&F^\prime\\
\downarrow&&\downarrow\\
B&=&B\end{array} $$
such that $U$ is $\phi$-twisted $T$-equivariant, i.e.
$U(p\cdot t)=U(p)\cdot\phi(t)$ for all $p\in F$, $t\in T$.\end{ddd}

Assume that $B$ is connected. We say that two $T$-principal bundles
over $B$ are twisted isomorphic if they are
$\phi$-twisted isomorphic for some (then uniquely determined) $\phi$.

\subsubsection{}
We consider a $T$-principal bundle $\pi:F\rightarrow B$.
We need the subgroups $S_0:=U(1)\times
\{1\}\subset T$ and $S_1:=\{1\}\times U(1)\subset T$.
We define $E_0:=F/S_0$ and $E_1:=F/S_1$. All these spaces fit into a diagram
\begin{equation}\label{cd1}\begin{array}{ccccc}
&&F&&\\
&p_0\swarrow&&p_1\searrow&\\
E_0&&\pi\downarrow&&E_1\\
&\pi_0 \searrow&& \pi_1 \swarrow\\
&&B&&
\end{array}
\ ,\end{equation}
where $p_i$ and $\pi_i$ are $U(1)$-principal bundles in a natural way.
We consider a class $h\in H^3(F,\Z)$.
\begin{ddd}
We say that the pair $(F,h)$ is dualizable, if
$h=p_0^*(h_0)+p_1^*(h_1)$ for some $h_i\in H^3(E_i,\Z)$.\end{ddd}

\subsubsection{}
We can now try to construct a $T$-dual of $(F,h)$ by iterated $T$-duality.
We first form the dual $({}_0 \hat F,{}_0\hat h)$ of the pair $(F\rightarrow
E_1,h)$. Note that we have the pull-back diagram
$$\begin{array}{rcr}F&\stackrel{p_0}{\rightarrow}&E_0\\
p_1=\pi_1^*\pi_0\downarrow&&\pi_0\downarrow\\
E_1&\stackrel{\pi_1}{\rightarrow}&B\end{array}\ .$$

Let $(\hat E_0,\hat h_0)$ be a dual of $(E_0,h_0)$.
Then we get ${}_0\hat F$ by the pull-back diagram
$$\begin{array}{rcr}{}_0\hat F&\stackrel{{}_0\hat p_0}{\rightarrow}&\hat E_0\\
{}_0\hat p_1=\pi_1^*\hat \pi_0\downarrow&&\hat \pi_0\downarrow\\
E_1&\stackrel{\pi_1}{\rightarrow}&B\end{array}\ .$$
Furthermore we get
$${}_0\hat
h:={}_0\hat p_0^*(\hat h_0)+{}_0 \hat p_1^*(h_1)\ .$$
Now we form the dual $(\hat F,\hat h)$ of the pair $({}_0\hat F\rightarrow
\hat E_0,{}_0\hat h)$.
Let $(\hat E_1,\hat h_1)$ be the dual of $(E_1,h_1)$.
Then we get  $\hat F$ by the pull-back

$$\begin{array}{rcr}\hat F&\stackrel{\hat p_0}{\rightarrow}&\hat E_0\\
\hat p_1\downarrow&&\hat \pi_0\downarrow\\
\hat E_1&\stackrel{\hat \pi_1}{\rightarrow}&B\end{array}\ .$$
 and
 \begin{equation}
\hat h=\hat p_0^*(\hat h_0)+\hat p_1(\hat h_1)\
.\label{eq:formula_for_h_hat}
\end{equation}

\subsubsection{}
Note that this construction of the iterated dual involves the choice of a representation
$h=p_0^*(h_0)+p_1^*(h_1)$. The goal of the following discussion is to
show that the bundle $\hat F\rightarrow B$ may depend on this choice
even if we consider it up to twisted equivalence. It should be
remarked that our examples with a non-unique dual do
not depend on the existence of torsion in cohomology and therefore
would also show up if we only worked with de Rham cohomology.

Our example should be contrasted with the constructions of \cite{bhm},
were a very similar definition of $T$-duality for torus bundles is
used, but in which case (at least according to the authors) the
$T$-dual (which exists under conditions similar to ours) is uniquely
determined (up to isomorphism).
In \cite{mr}, a different approach to $T$-duality for torus bundles
is used, based on continuous trace algebras over the initial bundle
and actions of $\reals^n$ on the continuous trace algebra. Under our
existence assumption, the construction of \cite{mr} also gives rise to
a classical dual torus bundle, which is claimed to be uniquely determined. The
relationship to our construction is not quite clear, we plan to
investigate this, and to give more information about the higher
dimensional case in a subsequent paper.

\subsubsection{}
A $T$-principal bundle $F\rightarrow B$ gives rise to Chern
classes $c_0,c_1\in H^2(B,\Z)$ of the bundles $E_0,E_1$. The pair
$(c_0,c_1)$ determines the isomorphism class of $F$ by the proof of
Corollary \ref{cor:twisted_iso_class}. 
We consider this pair of Chern classes as a class $c(F)\in
H^2(B,\Z^2)$ in the natural way. Then the Chern classes of
the dual are $-((\pi_0)_!(h_0),(\pi_1)_!(h_1))$.
Note that $GL(2,\Z)$ acts on the cohomology with coefficients in
$\Z^2$.

\subsubsection{}
Choose now $\phi\in GL(2,\Z)$. Then we can define a new $T$-principal
bundle ${}^\phi F$. It has the same underlying fiber bundle
$F\rightarrow B$, but we redefine the action of $T$ such that
$$\begin{array}{ccc}{}^\phi F&\stackrel{\id}{\rightarrow}&F^\prime\\
\downarrow&&\downarrow\\
B&=&B\end{array} \
,$$
is a $\phi$-twisted isomorphism.
Let $\sigma:GL(2,\Z)\rightarrow GL(2,\Z)$ be the bijection (of order two)
$$\left(\begin{array}{cc} a&b\\c&d\end{array}\right) \mapsto (ad-bc) \left(\begin{array}{cc} a&-c\\-b&d\end{array}\right)\ .$$
\begin{lem}\label{lem:how_c_changtes}
We have $c({}^\phi F)=\phi^\sigma c(F)$.
\end{lem}
\proof
Let $\phi=\left(\begin{array}{cc} a&b\\c&d\end{array}\right)$.
Then by an easy computation (which has only to be carried out for the
generators $
\left(\begin{smallmatrix}
  -1 & 0 \\ 0 & 1
\end{smallmatrix}\right)$, $\left(
\begin{smallmatrix}
  0 & -1\\ 1 & 0
\end{smallmatrix}\right)$, $\left(
\begin{smallmatrix}
  1 & -1\\ 1 & 0
\end{smallmatrix}\right)$
of $GL(2,\Z)$) we see that
${}^\phi E_0=(E_0^a\otimes E_1^{-c})^{\det(\phi)}$ and
${}^\phi E_1=(E_0^{-b}\otimes E_1^d)^{\det(\phi)}$.
Therefore $c({}^\phi F)=\det(\phi)\left(\begin{array}{cc}
    a&-c\\-b&d\end{array}\right)(c_0,c_1)$.
\hB

\begin{kor}\label{cor:twisted_iso_class}
The twisted isomorphism class of the $T$-principal bundle $F$ is
determined precisely by the orbit $GL(2,\Z)c(F)\subset H^2(B,\Z^2)$.
\end{kor}
\proof
  If $F$ and $F'$ are isomorphic $T$-principal bundles over $B$, then
  $c(F)=c(F')\in H^2(B,\Z^2)$. If $F$ and $F'$ are $\Phi$-twisted
  isomorphic, then $F$ and ${}^\Phi F'$ are isomorphic, and by Lemma
  \ref{lem:how_c_changtes} $c(F)\in Gl(2,\Z)c(F')$.

  Finally, if $F$ and $F'$ are two bundles with $c(F)=c(F')$, then the
  $U(1)$-principal bundles $E_0$ and $E_0'$ as well as $E_1$ and
  $E_1'$ are isomorphic. It follows that $F$, as the pullback of $E_1$
  along $E_0$, is isomorphic to the pullback of $E_1'$ along $E_0'$,
  i.e.~to $F'$. If $c(F)\in GL(2,\Z)c(F')$ then $c(F)=c({}^\Phi F')$
  for a suitable $\Phi\in GL(2,\Z)$ and therefore $F$ and ${}^\Phi F'$
  are isomorphic $T$-principal bundles, so that $F$ and $F'$ are
  twisted isomorphic.
\hB

\subsubsection{}
Let $L_i$ be the line bundles associated to $E_i$ and
$V:=L_0\oplus L_1$. Let further $s:S(V)\rightarrow B$ be the unit sphere
bundle. Then we have natural embeddings $i_0:E_0\rightarrow S(V)$ and
$i_1:E_1\rightarrow S(V)$. We have a decomposition
$$S(V)=D(L_0)\times_BE_1\cup E_0\times_B D(L_1)\ .$$
The associated Mayer-Vietoris sequence gives the exact sequence
\begin{equation}
H^2(F,\Z)\to H^3(S(V),\Z)\stackrel{i_0^*\oplus i_1^*}{\rightarrow} H^3(E_0,\Z)\oplus
H^3(E_1,\Z)\stackrel{p_0^*-p_1^*}{\rightarrow} H^3(F,\Z).
\label{eq:MV}
\end{equation}

Let $r_i\in H^3(E_i,\Z)$. Then we have  $p_0^*(r_0)+p_1^*(r_1)=0$ if
an only if $(r_0,-r_1)\in \im(i_0^*\oplus i_1^*)$. If this is satisfied we get a (second)
splitting of $0=p_0^*(0)=p_1^*(0)\in H^3(F,\Z)$.

To understand the corresponding dual, we compute $(\pi_{i})_{!}(r_i)$. The dual of
any $T$-bundle with splitting $0=p_0^*(0)+p_1^*(0)$ has Chern class
$(0,0)$. If we can find an example as above  with $((\pi_0)_!((i_0)^*X),(\pi_1)_!((i_1)^*X))\ne
0$, the latter can not lie in the $GL(2,\Z)$-orbit of $(0,0)$ and
therefore not even the underlying bundle of the second dual is twisted
isomorphic to the first one.

\subsubsection{}
Choose now $B=S^2$, $c_0=c_1$ the generator of $H^2(S^2,\Z)$. Then
$E_0=E_1$ has underlying space $S^3$ with the Hopf principal
fibration. 
In this case $(\pi_i)_!\colon H^3(S^3,\Z)\to H^2(S^2,\Z)$ is an
isomorphism. Moreover, $F$ is a $U(1)$-bundle over $S^3$ and
therefore $H^2(F,\Z)=0$, consequently $i_0^*\oplus i_1^*$ in
\eqref{eq:MV} is injective.
The Gysin sequence for $S(V)$ gives
$$H^3(B,\Z)\rightarrow H^3(S(V),\Z)\stackrel{s_!}{\rightarrow}
H^0(B,\Z)\stackrel{c_0\cup c_1}{\rightarrow} H^4(B,\Z)\ ,$$
i.e.~$H^3(S(V),\Z)\iso\Z\ne \{0\}$.

It follows that there is $0\ne X\in H^3(S(V),\Z)$ such that
$i_0^*(X)\oplus i_1^*(X)\ne 0$, and therefore
$(\pi_0)_!(i_0^*(X))\oplus (\pi_1)_!(i_1^*(X))\ne 0$ and we are done
by the above observation.

 \end{document}